\documentclass[11pt]{article}

\usepackage[mathscr]{eucal}
\usepackage{color}
\usepackage{amsmath}
\usepackage{amssymb}
\usepackage{amsfonts}
\usepackage{amscd}
\usepackage{amsthm}
\usepackage{latexsym}
\usepackage{graphicx}
\usepackage{subfig}

\def\be{\begin{equation}}
\def\ee{\end{equation}}
\def\bse{\begin{subequations}}
\def\ese{\end{subequations}}

\newtheorem{thm}{Theorem}

\newtheorem{lem}[thm]{Lemma}
\newtheorem{prop}[thm]{Proposition}

\newtheorem{rem}[thm]{Remark}

\def\bse{\begin{subequations}}
\def\ese{\end{subequations}}

\usepackage{geometry}

\title{Global existence and asymptotic behavior of solutions to a nonlocal Fisher-KPP type problem}

\author{Shen Bian\footnote{Beijing University of Chemical Technology, 100029, Beijing. Email: \texttt{bianshen66@163.com}. Partially supported by China Postdoctoral Science Foundation, No. 2014M560037.}
  \and Li Chen\footnote{Universit\"at Mannheim, 68131, Mannheim. Email: \texttt{chen@math.uni-mannheim.de}. Partially supported by the National
Natural Science Foundation of China (NSFC), No. 11271218.}
  \and Evangelos A. Latos\footnote{Universit\"at Mannheim, 68131, Mannheim. Email: \texttt{evangelos.latos@math.uni-mannheim.de}.}
}

\date{}

\begin{document}
\let\cleardoublepage\clearpage

\maketitle

\begin{abstract}
In this work, we consider a nonlocal Fisher-KPP reaction-diffusion
problem with Neumann boundary condition and nonnegative initial data
in a bounded domain in $\mathbb{R}^n (n \ge 1)$, with reaction
term $u^\alpha(1-m(t))$, where $m(t)$ is the total mass at time $t$. When $\alpha \ge 1$ and the initial mass
is greater than or equal to one, the problem has a unique
nonnegative classical solution. While if the initial mass is less
than one, then the problem admits a unique global solution for
$n=1,2$ with any $1 \le \alpha <2$ or $n \ge 3$ with any $1 \le
\alpha < 1+2/n$. Moreover, the asymptotic convergence to the
solution of the heat equation is proved. Finally, some numerical
simulations in dimensions $n=1,2$ are exhibited. Especially, for
$\alpha>2$ and the initial mass is less than one, our numerical
results show that the solution exists globally in time and the mass
tends to one as time goes to infinity.
\end{abstract}

\section{Introduction}

In this work we consider the following nonlocal initial boundary value problem,
\bse\label{nkpps0}
\begin{align}
&\ u_t-\Delta u=u^\alpha \left(1-\int_\Omega u(x,t)dx\right),\quad   &&x\in\Omega, t>0,
\\
&\ \nabla u\cdot \nu =0, &&x\in\partial \Omega,
\\
&\ u(x,0)=u_0(x)\geq0,\quad &&x\in\Omega,
\end{align}
\ese where $u$ is the density, $\Omega$ is a smooth bounded domain
in $\mathbb{R}^n$, $n\geq 1$, $\alpha\geq1$ and $\nu$ is the outer
unit normal vector on $\partial\Omega$. Without loss of generality,
throughout this paper we assume $|\Omega|=1$ (otherwise, rescale the
problem by $|\Omega|$), let $m(t)=\int_\Omega u(x,t)dx$ and
$m_0=m(0)$. A damping term with $\sigma>0$ can also be included to
get $u_t-\Delta u+\sigma u=u^\alpha(1-m(t))$. In this case,  similar
results to this paper can also be obtained. For simplicity, we
assume that $\sigma=0.$


In the 1930s, Fisher \cite{f37} and Kolmogorov, Petrovskii, Piskunov \cite{kpp} in population dynamics and Zeldovich, Frank-Kamenetskii \cite{zfk38} in combustion theory  started to study problems with this kind of reaction terms. Actually, they introduced the scalar reaction-diffusion equation
$
\frac{\partial u}{\partial t} = \frac{\partial^2u}{\partial x^2} +F(u),
$
 and studied the existence, stability and speed of propagation. In the theory of population dynamics, the function $F$ is considered as
the rate of the reproduction of the population. It is usually of the form
$$
F(u)=\beta u^\alpha(1-u)-\gamma u.
$$
From the above model two cases emerge depending on the values of $\alpha$.

In the case of $\alpha = 1$, the reproduction rate is proportional to the density $u$ of the population and to available resources $(1-u)$. The last term, $-\gamma u$, describes the mortality of the population.

The case $\alpha = 2$, which is the motivation for our work, considers the addition of sexual reproduction to the model with the reproduction rate proportional to the square of the density, see \cite{VVpp}. For more information on reaction-diffusion waves in biology, we refer to the review paper of Volpert and Petrovskii \cite{volpet}.

Next we pass to the relation between the local and the non-local consumption of resources. In the local reaction-diffusion problem
\be\label{lrdp}
\frac{\partial u}{\partial t} = \frac{\partial^2u}{\partial x^2} +\beta u^\alpha(1-u)-\gamma u,
\ee
where $u$ is the population density, $\frac{\partial^2u}{\partial x^2}$ describes the random displacement of the individuals of this population and the reaction term represents their reproduction and mortality. Moreover, the reaction term consists of  the reproduction term which is represented by the population density to a power, $u^\alpha$, multiplied with the term $(1-u)$ which stands for the local consumption of available resources.

The nonlocal version of the above problem is
\be\label{ide}
\frac{\partial u}{\partial t}=\frac{\partial^2u}{\partial x^2}+\beta u^\alpha\left(1-\int_{-\infty}^{\infty}\phi(x-y)u(y,t)dy\right)-\gamma u,
\ee
where $\beta,\gamma>0$  and $\int_{-\infty}^{\infty}\phi(y)dy=1.$ It can be seen as the case where the individual, located at a certain point, can consume resources in some area around that point. $\phi(x - y)$ represents the probability density function that describes the distribution of individuals around their average positions and it depends on the distance from the average point $x$ to the actual point $y$. One can easily verify that if $\phi$ is a Dirac $\delta$-function, then the nonlocal problem reduces to (\ref{lrdp}).

 In the current paper we will study problems with reaction terms similar to the above nonlocal reaction terms. There are some already known results on the reaction-diffusion equation with a nonlocal term,
$$u_t =\Delta u+F(t,u,I(u)),\quad I(u)=\int_\Omega u(y,t)dy,$$
in a bounded domain $\Omega$. However, compared to the local  version, the results for the nonlocal reaction terms of Fisher-KPP type are relatively limited. Here we list some of the known recent results.

 Anguiano, Kloeden and Lorenz considered $F = f(u)I(u^4)(1-I(u^4))$ and proved the existence of a global attractor \cite{akl10}. Wang and Wo \cite{ww11} proved the convergence to a stationary solution with $F = u^m - I(u^m),\ m > 1$. For $F  = \alpha e^{\gamma u} + bI(e^{\gamma u}),$ Pao \cite{Pao92} studied the existence or nonexistence of stationary solutions. Liu, Chen and Lu  \cite{lcl09} proved also the blow-up of solutions for a similar equation, see also \cite{dlx03}. Rouchon obtained global estimates of solutions in \cite{r03}. For more information on nonlocal KPP-Fisher type problems, we refer to a recent book by Volpert \cite{vol2}.

Nonlocal Fisher-KPP type reaction terms can describe also Darwinian evolution of a structured population density or the behavior of cancer cells with therapy as well as polychemotherapy and chemotherapy, we refer the interested reader to the models found in \cite{Lorz:2011hl,Lorz:2013vp,Lorz:2013hq}.

Bebernes and Bressan \cite{bb82} (see also Bebernes \cite{be82}, Pao\cite{Pao92}) considered the equation with reaction term
$$
F(t,u,I(u))=f (t,u(t,x)) + \int_\Omega g (t,u(t,y)) dy,\quad t>0,\quad x\in\Omega.
$$
They considered the case when $f(t,u) = e^u,\ g(t,u) = ke^u\ (k > 0)$, for which the above problem represents an ignition model for a compressible reactive gas, and proved that solutions blow-up.

Later, Wang and Wang \cite{ww96} considered a power-like nonlinearity, i.e.
$$
F(t,u,I(u))=\int_\Omega u^p(t,y)dy-ku^q(t,x),\quad t>0,\quad x\in\Omega,
$$
with $p, q > 1$, and proved the blow-up of the solutions.

 Budd, Dold and Stuart \cite{bds93}, Hu and Yin \cite{hy95} considered a similar to the above problem in the case $ p=2$ and general $p$ respectively,
$$
F(t,u,I(u))=u^p-\frac{1}{|\Omega|}\int_{\Omega}u^p(t,y)dy,\quad t>0,\quad x\in\Omega.
$$
With this typical structure, the energy of the solutions is
conserved (under Neumann boundary conditions). For this kind of
nonlocal problems it is known \cite{ww96} that there is no
comparison principle and they are the closest models to the ones we
are considering in this work. For a general study  on nonlocal
problems, we refer to Quittner's and Souplet's book \cite{sbook} as
well as the paper by Souplet \cite{s98}.

\subsection{Preliminary discussion}

In this article, we will focus on (\ref{nkpps0}) which has a
reaction term of the type
$$F(t,u,I(u))=u^\alpha\left(1-\int_\Omega u(x,t)dx\right)$$
for $\alpha\geq1$. One additional fact that makes this problem more difficult to handle is the lack of comparison principle as one can see for example from \cite{ww96}.

For nonnegative $u$, formally by integrating (\ref{nkpps0}) over $\Omega$, we get,
$$
m'(t)=(1-m(t))\int_\Omega u^\alpha dx,
$$
where $m(t)=\int_\Omega u(x,t)dx$ is the total mass at time $t$.

If we start at time $t_0$ such that $1-m(t_0)<0$, which means that $m(t_0)>1$, we can see that $m'(t)$ is negative and therefore $m(t)$ decreases in time. In this case it is natural to expect the global existence of solutions.

On the other hand, if we start at time $t_0$ such that $1-m(t_0)>0$,
we can see that  $m(t)$ increases in time. However if $(1-m(t))$
remains positive, the equation has a similar structure to the heat
equation with a power-like reaction term for which we know that
the problem might have no global solution for super-critical exponent $\alpha<1+2/n$ (see for
example \cite{bal1,bal2,fuj1,fuj2,fuj3,kapl}). In this paper we give
a negative answer to this observation. Our main
results are the following two theorems:

\begin{thm}\label{thm1}
Let $n \geq 1$, $\alpha \ge 1$ and $\int_{\Omega} u_0(x) dx=m_0>0$. Assume $u_0$ is nonnegative and $u_0 \in L^k(\Omega)$ for any $1<k<\infty$. Then for $m_0<1$ with $\alpha$ satisfying
   \begin{align}
    1 \le \alpha<1+2/n, ~~&n \ge 3, \\
    1 \le \alpha<2, ~~&n=1,2,
   \end{align}
or $m_0 \ge 1$ with arbitrary $\alpha \ge 1$, problem (\ref{nkpps0})
has a unique nonnegative classical solution. Moreover, the following
a priori estimates hold true. That's for $m_0<1$,
    \begin{align}
    \|u(\cdot,t)\|_{L^k(\Omega)} \le
    C+C~t^{-\frac{k-1}{\alpha-1}}~~\mbox{ for
    any} ~t>0.
    \end{align}
For $m_0 \ge 1$,
   \begin{align}
    \|u\|_{L^k(\Omega)}^k \le C+ C~ t^{-\frac{n(k-1)}{2}},\quad n \ge 3, \\
    \|u\|_{L^k(\Omega)}^k \le C+ C ~ t^{-(k-1)},\quad n=1,2.
    \end{align}
Here $C$ denote different constants depending on $m_0,k,\alpha$, but not depending on $\|u_0\|_{L^k(\Omega)}$.
\end{thm}
\begin{thm}\label{thm2}
Let $u(x,t)$ be the unique nonnegative classical solution obtained
from Theorem \ref{thm1}, $v$ be the solution to the heat equation
with Neumann boundary condition and initial data $\int_\Omega
v_0(x)dx= m_0$, then as $t \to \infty$,
\begin{align}
\|u(\cdot,t)-v(\cdot,t)-(m_0-1)\|_{L^2(\Omega)}\leq C_1e^{-C_2t},
\end{align}
where $C_1,C_2$ are constants depending on the initial mass $m_0$ and $\|u_0\|_{L^{2 \alpha}(\Omega)}$.
\end{thm}

This paper is organized as follows. In Section 2, we firstly present
the dynamics of the mass. The global existence of the solutions and thus the proof the Theorem \ref{thm1} are shown in Section 3.
Section 4 is devoted to the proof of Theorem \ref{thm2}.
Section 5 shows the numerical simulations of the problem which give
the motivation for our future study for the case $\alpha \ge 1+2/n$.
Finally, Section 6 concludes the main work of our paper, some open
questions for problem (\ref{nkpps0}) are also addressed.

\section{Dynamics of the total mass}
The evolution of the total mass plays the key role in our proof for the global existence of the classical solution. Therefore, we firstly give the evolution of mass
$\int_{\Omega}u(t)dx$ in time.
\begin{lem}\label{massestimate}
For $m_0>0$, the mass $\int_{\Omega} u(t) dx=m(t)$ satisfies
\begin{align} \label{massbdd}
\min\{ 1,m_0  \} \le m(t) \le \max \{ 1,m_0 \}.
\end{align}
Furthermore, we have the following decay estimates
\begin{align}\label{massdecay}
|1-m(t)|\leq |1-m_0|e^{-\min\{1,m_0^\alpha\}t}.
\end{align}
\end{lem}
\noindent {\textbf{Proof.}} We return to the original problem
(\ref{nkpps0}) and integrate it over $\Omega$ to get: \be\label{l1s0}
m'(t)=(1-m(t))\int_\Omega u^\alpha dx. \ee There are two
possibilities depending on the initial mass.
\begin{itemize}
\item If we start at time $t_0$ where $1-m(t_0)>0$, we can see that $m'(t)$ is positive and therefore $m$ increases in time. Moreover, with the use of Jensen's inequality we get by using $m(t)\geq m_0,$
$$m'(t)\geq (1-m(t))m^\alpha(t)\geq
(1-m(t))m_0^{\alpha}.
$$
By solving this inequality we get a lower bound on the speed with
which $m(t)$ increases to 1 i.e.
$$
m(t)\geq1-e^{-m_0^\alpha t}(1-m_0).
$$
\item If $1-m(t_0)<0$,  then $m$  decreases in time. By monotonicity, we get  $m(t)\leq m_0$ and again by Jensen's inequality,
$$
m'(t)\leq(1-m(t))m^\alpha(t)<1-m(t),
$$
then we get that
$$
m(t)\leq1+(m_0-1)e^{-t}.
$$
\end{itemize}
By putting together the above two cases we have the expected results. $\Box$

\section{Global Existence}
This section mainly focuses on the global existence of the classical
solution to (\ref{nkpps0}). We will use the following
ODE inequality from \cite{BL13}, which was also used in \cite{BL14}.

\begin{lem}\label{BL14ode}
Assume $y(t) \ge 0$ is a $C^{1}$ function for $t>0$ satisfying
$$y'(t)\le \alpha- \beta y(t)^a$$ for $a>1,~\alpha > 0,~\beta>0 $, then
$y(t)$ has the following hyper-contractive property
\begin{align} \label{agle1}
    y(t) \le (\alpha/\beta)^{1/a}+ \left[\frac{1}{\beta(a-1)t}\right]^{\frac{1}{a-1}}  \quad \mbox{ for any ~} t >
    0.
\end{align}
Furthermore, if $y(0)$ is bounded, then
\begin{align} \label{ytlessy0}
 y(t) \le \max \left( y(0),  (\alpha/\beta)^{1/a} \right).
\end{align}
\end{lem}

The proof of global existence heavily depends on the following two
\emph{a priori} estimates, Proposition \ref{globalsub} and
Proposition \ref{globalsuper}, and then we will use the compactness
arguments to close the proof. Due to the preliminary discussion, we
will divide the a priori estimates into two cases $m_0 < 1$ and $m_0
\ge 1$. Firstly, we focus on $m_0 < 1$.

\begin{prop}\label{globalsub}
Let $n \ge 1$ and $m_0 < 1$. If $\alpha$ satisfies
   \begin{align*}
    1 \le \alpha<1+2/n, ~~&n \ge 3, \\
    1 \le \alpha<2, ~~&n=1,2,
   \end{align*}
   then for any $1<k<\infty$, the nonnegative solution of (\ref{nkpps0}) satisfies
    \begin{align}
    \|u(\cdot,t)\|_{L^k(\Omega)}\le
    C(m_0,k,\alpha)+C(m_0,k,\alpha)~t^{-\frac{k-1}{\alpha-1}}~~\mbox{for any}~t>0.
    \end{align}
    Moreover, if $u_0(x) \in L^k(\Omega)$, then
    \begin{align}
    \|u(\cdot,t)\|_{L^k(\Omega)} \le \max\left\{
    \|u_0(x)\|_{L^k(\Omega)},C(m_0,k,\alpha) \right\},
    \end{align}
    and for any $0<T<\infty$
    \begin{align*}
        \nabla u^{\frac{k}{2}} \in L^2 \left(0,T;L^2(\Omega)\right).
    \end{align*}
\end{prop}
\noindent\textbf{ Proof.} Since $m_0 < 1$, by lemma
\ref{massestimate} one has $m_0 \le m(t) \le 1$. Using $ku^{k-1}$ as a test
function for equation (\ref{nkpps0}) and integrating it by parts
\begin{align}
\frac{d}{dt} \int_{\Omega}  u^k dx = - \frac{4(k-1)}{k}
\int_{\Omega} |\nabla u^{\frac{k}{2}} |^2 dx+ k \int_{\Omega}
u^{k+\alpha-1} dx \left(1-\int_{\Omega} u dx \right), \nonumber\\
\frac{d}{dt} \int_{\Omega}  u^k dx + \frac{4(k-1)}{k} \int_{\Omega}
|\nabla u^{\frac{k}{2}} |^2 dx+ k m(t) \int_{\Omega} u^{k+\alpha-1}
dx = k \int_{\Omega} u^{k+\alpha-1} dx. \label{guji2}
\end{align}
Choosing $1<k'<k+\alpha-1,$ combining H\"{o}lder's inequality and the Sobolev embedding theorem
one has
\begin{align}
&\int_{\Omega} u^{k+\alpha-1} dx  = \int_{\Omega}  u^{\lambda \frac{k}{2}\frac{2(k+\alpha-1)}{k}} u^{(1-\lambda) \frac{k}{2}\frac{2(k+\alpha-1)}{k}}  dx \nonumber\\
\leq &\left\|  u^{\frac{k}{2}}
\right\|_{L^{p}(\Omega)}^{\lambda\frac{2(k+\alpha-1)}{k}}  \left \|   u^{\frac{k}{2}}
\right\|_{L^\frac{2k'}{k}(\Omega)}^{(1-\lambda)\frac{2(k+\alpha-1)}{k}} \nonumber \\
\leq & C(k)  \left( \left\|  \nabla u^{\frac{k}{2}}
\right\|_{L^2(\Omega)}^{\lambda}  \left \|   u^{\frac{k}{2}}
\right\|_{L^\frac{2k'}{k}(\Omega)}^{1-\lambda} + \left\|  u^{\frac{k}{2}}
\right\|_{L^2(\Omega)}^{\lambda}  \left \|   u^{\frac{k}{2}}
\right\|_{L^\frac{2k'}{k}(\Omega)}^{1-\lambda}
\right)^{\frac{2(k+\alpha-1)}{k}} \nonumber\\
\leq & C(k) \left( \left \|  \nabla u^{\frac{k}{2}}
\right\|_{L^2(\Omega)}^{\frac{2\lambda(k+\alpha-1)}{k}}  \left\|
u^{\frac{k}{2}}
\right\|_{L^\frac{2k'}{k}(\Omega)}^{\frac{2(1-\lambda)(k+\alpha-1)}{k}}
+ \left\| u^{\frac{k}{2}}
\right\|_{L^2(\Omega)}^\frac{2\lambda(k+\alpha-1)}{k}  \left\| u^{\frac{k}{2}}
\right\|_{L^\frac{2k'}{k}(\Omega)}^\frac{2(1-\lambda)(k+\alpha-1)}{k} \right),
\label{guji3}
\end{align}
where $\lambda$ is the exponent from H\"older's inequality, i.e.
\begin{align}
\lambda=\frac{ \frac{k}{2k'}- \frac{k}{2(k+\alpha-1)} }{
\frac{k}{2k'} -\frac{1}{p} } \in (0,1),
\end{align}
and $p$ satisfies
\begin{align}
\left\{
    \begin{array}{ll}
      p=\frac{2n}{n-2}, & n \ge 3, \\
     \frac{2(k+\alpha-1)}{k}<p<\infty , & n=2, \\
      p=\infty, & n=1.
    \end{array}
  \right.
\end{align}
Now we will divide the analysis into three cases $n\ge 3,n=2$ and $n=1$.

For $n \ge 3$, $p=\frac{2n}{n-2}$ and then
\begin{align}\label{lambda}
\lambda=\frac{ \frac{kn}{2k'}- \frac{kn}{2(k+\alpha-1)} }{
\frac{kn}{2k'} +1-\frac{n}{2} } \in (0,1),
\end{align}
with $k>\max \left\{ \frac{(n-2)(\alpha-1)}{2},1\right\}$. Taking
$k'>\frac{(\alpha-1) n}{2}$, simple computations arrive at
\begin{align*}
\frac{2 \lambda (k+\alpha-1)}{k} & = \frac{\frac{kn}{k'}+
\frac{(\alpha-1) n}{k'}-n  }{ \frac{kn}{2 k'}+ 1-\frac{n}{2} }  < 2.
\end{align*}
To sum up, for $k'>\max\left\{ \frac{(\alpha-1) n}{2},1 \right\},$
thanks to the Young's inequality, from (\ref{guji3}) one has
\begin{align}
\int_{\Omega} u^{k+\alpha-1} dx \le &\frac{k-1}{k^2} \left\|\nabla
u^{\frac{k}{2}} \right\|_{L^2(\Omega)}^2 + C(k) \left\|
u^{\frac{k}{2}} \right\|_{L^\frac{2k'}{k}(\Omega)}^{(1-\lambda)
\frac{2(k+\alpha-1)}{k}
\frac{1}{1-\frac{\lambda(k+\alpha-1)}{k}}} \nonumber\\
& +C(k) \left\| u^{\frac{k}{2}}
\right\|_{L^2(\Omega)}^\frac{2\lambda(k+\alpha-1)}{k}  \left\| u^{\frac{k}{2}}
\right\|_{L^\frac{2k'}{k}(\Omega)}^\frac{2(1-\lambda)(k+\alpha-1)}{k}.
\label{guji4}
\end{align}
Letting
\begin{align}
r=(1-\lambda) \frac{2(k+\alpha-1)}{k}
\frac{1}{1-\frac{\lambda(k+\alpha-1)}{k}},
\end{align}
recalling $m(t) \ge m_0$, together (\ref{guji2}) with (\ref{guji4}) we arrive at
\begin{align}
& \frac{d}{dt} \int_{\Omega}  u^k dx + k m_0 \int_{\Omega}
u^{k+\alpha-1} dx + \frac{3(k-1)}{k} \| \nabla u^{\frac{k}{2}} \|_{L^{2}(\Omega)}^2 \nonumber\\
\le ~ & C(k) \|u\|_{L^{k'}(\Omega)}^{\frac{kr}{2}}+C(k) \|u\|_{L^{k}(\Omega)}^{\lambda(k+\alpha-1)}
\|u\|_{L^{k'}(\Omega)}^{(1-\lambda)(k+\alpha-1)}. \label{guji5}
\end{align}
On the other hand, using H\"{o}lder's inequality with
$1<k'<k+\alpha-1$ we have
\begin{align}
\|u\|_{L^{k'}(\Omega)} \le C
\|u\|_{L^{k+\alpha-1}(\Omega)}^{\theta}\|u\|_{L^1(\Omega)}^{1-\theta},  \\
\|u\|_{L^{k}(\Omega)} \le C
\|u\|_{L^{k+\alpha-1}(\Omega)}^{\eta}\|u\|_{L^1(\Omega)}^{1-\eta},
\end{align}
where
\begin{align}
\theta=\frac{(k+\alpha-1)(k'-1)}{k'(k+\alpha-2)} \in (0,1),~~\eta=\frac{(k+\alpha-1)(k-1)}{k(k+\alpha-2)} \in (0,1).
\end{align}
Hence
\begin{align}
\|u\|_{L^{k'}(\Omega)}^{\frac{kr}{2}} & \le \left(C
\|u\|_{L^{k+\alpha-1}(\Omega)}^{\theta}\|u\|_{L^1(\Omega)}^{1-\theta}
\right)^{\frac{kr}{2}} \le C(m_0,k) \|u\|_{L^{k+\alpha-1}(\Omega)}^{\frac{kr
\theta}{2}}. \label{inter}
\end{align}
Taking (\ref{guji5})-(\ref{inter}) into account we obtain that
\begin{align}
&\frac{d}{dt} \int_{\Omega}  u^k dx + k m_0 \|u\|_{L^{k+\alpha-1}(\Omega)}^{k+\alpha-1} + \frac{3(k-1)}{k} \| \nabla u^{\frac{k}{2}} \|_{L^{2}(\Omega)}^2 \nonumber\\
\le ~ & C(m_0,k) \|u\|_{L^{k+\alpha-1}(\Omega)}^{\frac{kr
\theta}{2}}+C(m_0,k) \|u\|_{L^{k+\alpha-1}(\Omega)}^{(k+\alpha-1)
[\lambda \eta + (1-\lambda) \theta]}. \label{guji6}
\end{align}
Here
\begin{align*}
\frac{kr}{2} & =  (1-\lambda) (k+\alpha-1)
 \frac{1}{1-\frac{\lambda(k+\alpha-1)}{k}}, \\
\frac{k+\alpha-1}{\theta} & =\frac{k+\alpha-2}{1-\frac{1}{k'}}.
\end{align*}
Recalling the definition of $\theta,\eta,\lambda$, direct computations show that
$ \lambda \eta + (1-\lambda) \theta <1$.
For
\begin{align}
1 \le \alpha < 1+\frac{2}{n},
\end{align}
one can derive that
\begin{align}\label{xiaoyu}
\frac{kr \theta}{2}< k+\alpha-1.
\end{align}
Next for $n=2$, $\frac{2(k+\alpha-1)}{k}<p<\infty$, by proceeding the similar arguments to the case $n \ge 3$ from (\ref{lambda}) to (\ref{guji6}), we obtain that for
\begin{align}
1 \le \alpha< 2-\frac{2}{p},
\end{align}
(\ref{xiaoyu}) holds true. When $n=1$, $p=\infty$, then for $1 \le \alpha<2$, (\ref{xiaoyu}) also holds true.

Therefore, combining the three cases $n \ge 3$, $n=2$ and $n=1$, using Young's inequality we obtain from
(\ref{guji6}) that
\begin{align}\label{Lk}
&\frac{d}{dt} \int_{\Omega}  u^k dx + k m_0 \|u\|_{L^{k+\alpha-1}(\Omega)}^{k+\alpha-1} +\frac{3(k-1)}{k} \| \nabla u^{\frac{k}{2}} \|_{L^{2}(\Omega)}^2   \nonumber\\
\le ~ & \frac{k m_0}{4}
\|u\|_{L^{k+\alpha-1}(\Omega)}^{k+\alpha-1} +\frac{k m_0}{4}
\|u\|_{L^{k+\alpha-1}(\Omega)}^{k+\alpha-1} +C(m_0,k).
\end{align}
In addition, H\"{o}lder's inequality yields that
\begin{align}
\left(\|u\|_{L^k(\Omega)}^{k}\right)^{1+\frac{\alpha-1}{k-1}}  \le
\|u\|_{L^{k+\alpha-1}(\Omega)}^{k+\alpha-1} \|u\|_{L^1(\Omega)}^{\frac{\alpha-1}{k-1}}.
\end{align}
Then using lemma \ref{BL14ode}, we
solve the ODE inequality
\begin{align} \label{ODEinequality}
\frac{d}{dt} \int_{\Omega}  u^k dx + \frac{k m_0}{2
\|u\|_{L^1(\Omega)}^{\frac{\alpha-1}{k-1}}} \left( \int_{\Omega} u^k
dx \right)^{1+\frac{\alpha-1}{k-1}} \le C(m_0,k)
\end{align}
to obtain that for any $1<k<\infty$,
\begin{align}\label{fartime}
\|u\|_{L^k(\Omega)}^k \le C(m_0,k,\alpha)+ \left[
\frac{C(m_0,k,\alpha)}{t} \right]^{\frac{k-1}{\alpha-1}}~~\mbox{ for
any} ~t>0.
\end{align}
Furthermore, if $u_0(x) \in L^k(\Omega)$ for any $1<k<\infty$, then taking
$y(t)=\int_{\Omega} u^k dx$ in lemma \ref{BL14ode} one has
\begin{align}\label{neartime}
\|u\|_{L^k(\Omega)} \le \max \left\{
\|u_0(x)\|_{L^k(\Omega)},C(m_0,k,\alpha) \right\}.
\end{align}
Now we integrate (\ref{Lk}) from $0$ to $T$ in time, then we can obtain that for any $T>0$
\begin{align*}
\int_{\Omega} u^k(T) dx+\int_0^{T}\int_{\Omega} \left| \nabla u^{\frac{k}{2}} \right|^2 dx dt + \int_0^{T}\int_{\Omega} u^{k+\alpha-1} dx dt   \le
\int_{\Omega} u_0^k dx+ C(m_0,k)T,
\end{align*}
from which we derive
\begin{align*}
\nabla u^{\frac{k}{2}} \in L^2 \left(0,T;L^2(\Omega)\right) \quad \mbox{for any} ~T>0.
\end{align*}
This completes the proof. $\Box$

For $m_0 \ge 1$, owing to lemma \ref{massestimate}, we know $m(t) \ge 1$
for any $t>0$, thus we have the following result
\begin{prop}\label{globalsuper}
Let $n \ge 1$ and $\alpha \ge 1$. Assume $u_0 \in L_+^1(\Omega)$ and $\int_{\Omega}u_0(x) dx=m_0 \ge 1$,
Then for any $1<k<\infty$, the nonnegative solution of (\ref{nkpps0}) satisfies that for any $t>0$
\begin{align}
\|u\|_{L^k(\Omega)}^k \le C(m_0,k)+ \left[
\frac{C(m_0,k)}{t} \right]^{\frac{n(k-1)}{2}},\quad n \ge 3, \label{38} \\
\|u\|_{L^k(\Omega)}^k \le C(m_0,k)+ \left[
\frac{C(m_0,k)}{t} \right]^{k-1},\quad n=1,2. \label{39}
\end{align}
Moreover, if $u_0 \in L^k(\Omega)$, then
\begin{align}
\int_0^{\infty}\int_{\Omega} \left| \nabla u^{\frac{k}{2}} \right|^2 dx dt  \le
\int_{\Omega} u_0^k dx.
\end{align}
\end{prop}
\noindent{\textbf{ Proof.}} Recalling lemma \ref{massestimate}, we know
that if $m_0>1$, then $m(t) \ge 1$ for any $t>0$. Hence the $L^k$
estimates (\ref{guji2}) can be reduced to
\begin{align}\label{k2}
\frac{d}{dt} \int_{\Omega}  u^k dx + \frac{4(k-1)}{k} \int_{\Omega}
|\nabla u^{\frac{k}{2}} |^2 dx \le 0.
\end{align}
For $n \ge 1$, using H\"{o}lder's inequality one has that for any $1<k<\infty$, the following estimate holds
\begin{align}
\|u\|_{L^k(\Omega)}^k & \le \|u\|_{L^{\frac{kp}{2}}(\Omega)}^{k\theta} \|u\|_{L^1(\Omega)}^{k(1-\theta)} \nonumber\\
& = \|u^{\frac{k}{2}}\|_{L^p(\Omega)}^{2 \theta} \|u\|_{L^1(\Omega)}^{k(1-\theta)}, \label{ktheta}
\end{align}
where $\theta=\frac{k-1}{k-\frac{2}{p}}$ and
\begin{align}
\left\{
    \begin{array}{ll}
      p=\frac{2n}{n-2}, & n \ge 3, \\
     2<p<\infty , & n=2, \\
      p=\infty, & n=1.
    \end{array}
  \right.
\end{align}
Thanks to the Sobolev embedding theorem and Young's inequality, from (\ref{ktheta}) one has
\begin{align}
\left(\|u\|_{L^k(\Omega)}^k \right)^{\frac{1}{\theta}} & \le C(n) \left( \|\nabla u^{\frac{k}{2}} \|_{L^2(\Omega)}^2 + \| u^{\frac{k}{2}} \|_{L^2(\Omega)}^2  \right) \|u\|_{L^1(\Omega)}^{\frac{k(1-\theta)}{\theta}}   \nonumber \\
& \le C(m_0,n) \|\nabla u^{\frac{k}{2}} \|_{L^2(\Omega)}^2 + C(m_0,n) \|u\|_{L^k(\Omega)}^k \nonumber \\
& \le C(m_0,n) \|\nabla u^{\frac{k}{2}} \|_{L^2(\Omega)}^2+ \frac{1}{2} \left(\|u\|_{L^k(\Omega)}^k \right)^{\frac{1}{\theta}} +C(m_0,n,k).
\end{align}
Plugging the above estimates into (\ref{k2}) yields that
\begin{align}\label{k3}
\frac{d}{dt} \int_{\Omega}  u^k dx +C(m_0,n,k) \left( \int_{\Omega} u^k
dx \right)^{\frac{1}{\theta}}+ C(m_0,n,k) \|\nabla u^{\frac{k}{2}} \|_{L^2(\Omega)}^2 \le C(m_0,k,n),
\end{align}
solving the ODE inequality we have that for any $t>0$
\begin{align}
\|u\|_{L^k(\Omega)}^k \le C(m_0,k,n)+ \left[
\frac{C(m_0,k,n)}{t} \right]^{\frac{n(k-1)}{2}},\quad n \ge 3, \\
\|u\|_{L^k(\Omega)}^k \le C(m_0,k,n)+ \left[
\frac{C(m_0,k,n)}{t} \right]^{k-1},\quad n=1,2.
\end{align}
Moreover, if $u_0 \in L^k(\Omega)$, then (\ref{k2}) directly yields
\begin{align}
\|u\|_{L^{k}(\Omega)}\le \|u_0\|_{L^k(\Omega)}.
\end{align}
Next integrating (\ref{k2}) from $0$ to $\infty$ in time we obtain that
\begin{align}
\int_0^{\infty}\int_{\Omega} \left| \nabla u^{\frac{k}{2}} \right|^2 dx dt  \le
\int_{\Omega} u_0^k dx.
\end{align}
This closes the proof. $\Box$
\begin{rem}
In fact, (\ref{38}) and (\ref{39}) also hold true for heat equation, and the uniform boundedness in time of the $L^k$ norm depends only on the initial mass, not depends on the initial $L^k$ norm.
\end{rem}
Now we have obtained the necessary a priori estimates to complete
the proof of  Theorem \ref{thm1}. It can be proved by standard
methods and for the convenience of the reader we mention the key
steps in the following. First of all, from the above estimates, we
can take $k=2$ and $k=2\alpha$ in Proposition \ref{globalsub} and
\ref{globalsuper} to get the estimates for $\|\nabla
u\|_{L^2(L^2(0,T))}$ and $\|u_t\|_{L^2(H^{-1}(0,T))}$ for any $T>0$.
By Aubin-Lions lemma \cite{SimJac86,XiuAnsgar} , we have the strong
compactness of $u$ in $L^2$ so that the nonlinear terms can be
handled. Therefore, the global existence of weak solutions (in the
sense of distributions) can be obtained by standard compactness
argument. Secondly, from the estimates of the weak solution in
Proposition \ref{globalsub} and Proposition \ref{globalsuper}, the
nonlinear term $u^{\alpha}(1-m(t))\in L^k([0,T]\times\Omega)$,
$\forall k>1$ for any $T>0$. The solution is a strong $W^{2,1}_k$
solution from classical parabolic theory, \cite{lsv,lb96}. By
Sobolev embedding, we can bootstrap it to get that classical
solution. In the end, the uniqueness can be obtained directly from
comparison principle \cite{lsv,lb96} since $u^{\alpha-1}(1-m(t))$ is
bounded from below.  $\Box$

\section{The long time behavior of solutions}
As we can see from the above arguments, equation (\ref{nkpps0})
has a unique classical solution. In this section, we will detect
the long time behavior of the global solution.
\begin{thm}\label{thmlong}
Assume $u_0 \in L^{k}(\Omega)$ for any $1<k<\infty$. Let $u$ be the classical solution to problem (\ref{nkpps0}) and $v$
be the solution to the heat equation with Neumann boundary
condition and initial data $v_0$ such that $\int_\Omega v_0(x)dx=
m_0$, Then as $t \to \infty$
$$
\|u(\cdot,t)-v(\cdot,t)-(m_0-1)\|_{L^2(\Omega)}\leq C_1e^{-C_2t},
$$
where the constants $C_1$, $C_2$ depend on $m_0$ and
$\|u_0\|_{L^{2\alpha}(\Omega)}$.
\end{thm}
\indent\textbf{Proof.} The difference between the two equations is
$$
(u-v)_t+\Delta(u-v)=u^\alpha(1-m(t)).
$$
Let $\overline u (t) =\int_\Omega u(x,t)dx$ and $\overline v (t) =\int_\Omega v(x,t)dx$.
By (\ref{l1s0}), we have $\overline{u}_t=m'(t)=(1-m(t))\int_\Omega u^\alpha dx$. $\overline{v}_t(t)=0$ because of $\overline{v}(t)=\overline{v}_0$. Therefore,
$$
(u-v)_t-(\overline{u}-\overline{v})_t+\Delta(u-v)=u^\alpha(1-m(t))-(1-m(t))\int_\Omega u^\alpha.
$$
The standard $L^2$ estimate shows that,
$$
\frac12\frac{d}{dt}\int_\Omega|(u-v)-(\overline{u}
-
\overline{v})|^2dx+\int_\Omega|\nabla (u-v)|^2dx
=
(1-m(t))\int_\Omega
\left[
\Big(
u^\alpha-\int_\Omega u^\alpha dy
\Big)
\big((u-v)-(\overline{u}-\overline{v})\big)
\right]
dx.
$$
By taking $k=2\alpha$ in Proposition \ref{globalsub} and Proposition
\ref{globalsuper}, we get
\begin{eqnarray*}
&&
\frac12\frac{d}{dt}\int_\Omega|(u-v)-(\overline{u}
-
\overline{v})|^2dx+\int_\Omega|\nabla (u-v)|^2dx\\
&\leq &
\frac12|1-m(t)|\int_\Omega|(u-v)-(\overline{u}-\overline{v})|^2dx
+ C |1-m(t)|,
\end{eqnarray*}
where $C$ depend on $m_0,\alpha$ and
$\|u_0\|_{L^{2\alpha}(\Omega)}$. Applying Poincar\'{e} inequality
and lemma \ref{massestimate}, we get
\begin{eqnarray*}
&&\frac12\frac{d}{dt}\int_\Omega|(u-v)-(\overline{u}
-\overline{v})|^2dx
+C(\Omega)
\int_\Omega|(u-v)-(\overline{u}-\overline{v})|^2 dx\\
&\leq & Ce^{-Ct}\int_\Omega|(u-v)-(\overline{u} -\overline{v})|^2dx+
C e^{-Ct}.
\end{eqnarray*}
From the above ODE we get the following estimate,
\begin{align*}
\int_\Omega|(u-v)-(\overline{u}-\overline{v})|^2 dx \le C_1 e^{-C_2 t},
\end{align*}
where $C_1,C_2$ are constants depending on $m_0$ and $\|u_0\|_{L^{2\alpha}(\Omega)}$.
Thus completes Theorem \ref{thmlong}. $\Box$



\section{Numerical Results}

 For $1 \le \alpha<2 ~(n=1,2)$, if $m_0<1$, it has been shown in the previous
sections that the solution will exist globally without any
restriction on the initial data. An interesting question is whether
the solution also exists globally for $\alpha>2$. In the following, we will give a complete numerical study for $n=1,2$ with any $\alpha \ge 1$.

\subsection{Numerical scheme}

For the numerical simulation, we consider the 2-dimensional equation
in $\Omega=[0,b]^2$ with any $b>0$
\begin{align}\label{semilinear}
\left\{
      \begin{array}{ll}
      u_t = \Delta u + u^{\alpha} \left(1-\iint_{\Omega} udxdy \right), ~~ (x,y) \in [0,b]\times[0,b],~t\geq 0, \\
      u(t=0) =u_{0}(x,y) \ge 0, ~~ (x,y) \in [0,b]\times[0,b],   \\[2mm]
      \frac{\partial u}{\partial x}(0,y,t)= \frac{\partial u}{\partial x}(b,y,t) = 0,
\\[2mm]
      \frac{\partial u}{\partial y}(x,0,t)= \frac{\partial u}{\partial y}(x,b,t) =
      0.
      \end{array}\right.
\end{align}
Denote
\begin{align} \label{fu}
f(u)= u^\alpha \left(1-\iint_{\Omega} u dxdy \right).
\end{align}
Here an alternating direction implicit difference scheme is applied
to construct the numerical computations.

Let $h$ be the space step and $\tau$ be the time step, $N=b/h+1$ is
the number of the discrete points, $T$ is the final time and
$K=T/\tau+1$.  Denote
\begin{align*}
    &x(i)=(i-1)*h, ~~i=1,2,3,\cdots,N-1,N, \\
    &y(j)=(j-1)*h, ~~j=1,2,3,\cdots,N-1,N, \\
    &t(k)=(k-1)*\tau,~~k=1,2,\cdots,T/\tau+1,\\
    &\Omega_h=\Big\{ (x(i),y(j)) | 1 \le i,j \le N \Big\}.
\end{align*}
The discrete solution at each time is presented as a matrix
$u_{i,j}^k \in \mathbb{R}^{N \times N}$, where
\begin{align*}
    &u_{i,j}=u\Big(x(i),y(j),\cdot\Big), ~~  1 \le i,j \le N,  \\
    &u^k=u\Big(\cdot,\cdot,t(k)\Big),~~k=1,2,\cdots,K.
\end{align*}
Next we introduce some notations:
\begin{align*}
    \delta_x^2 u_{ij}=\frac{u_{i+1,j}-2 u_{ij}+u_{i-1,j}}{h^2},\\
    \delta_y^2 u_{ij}=\frac{u_{i,j+1}-2 u_{ij}+u_{i,j-1}}{h^2}.
\end{align*}

\begin{figure}[htbp] \centering
 \subfloat[$u(x,y)$ at t=0.5s] {
     \includegraphics[height=7cm,width=7cm]{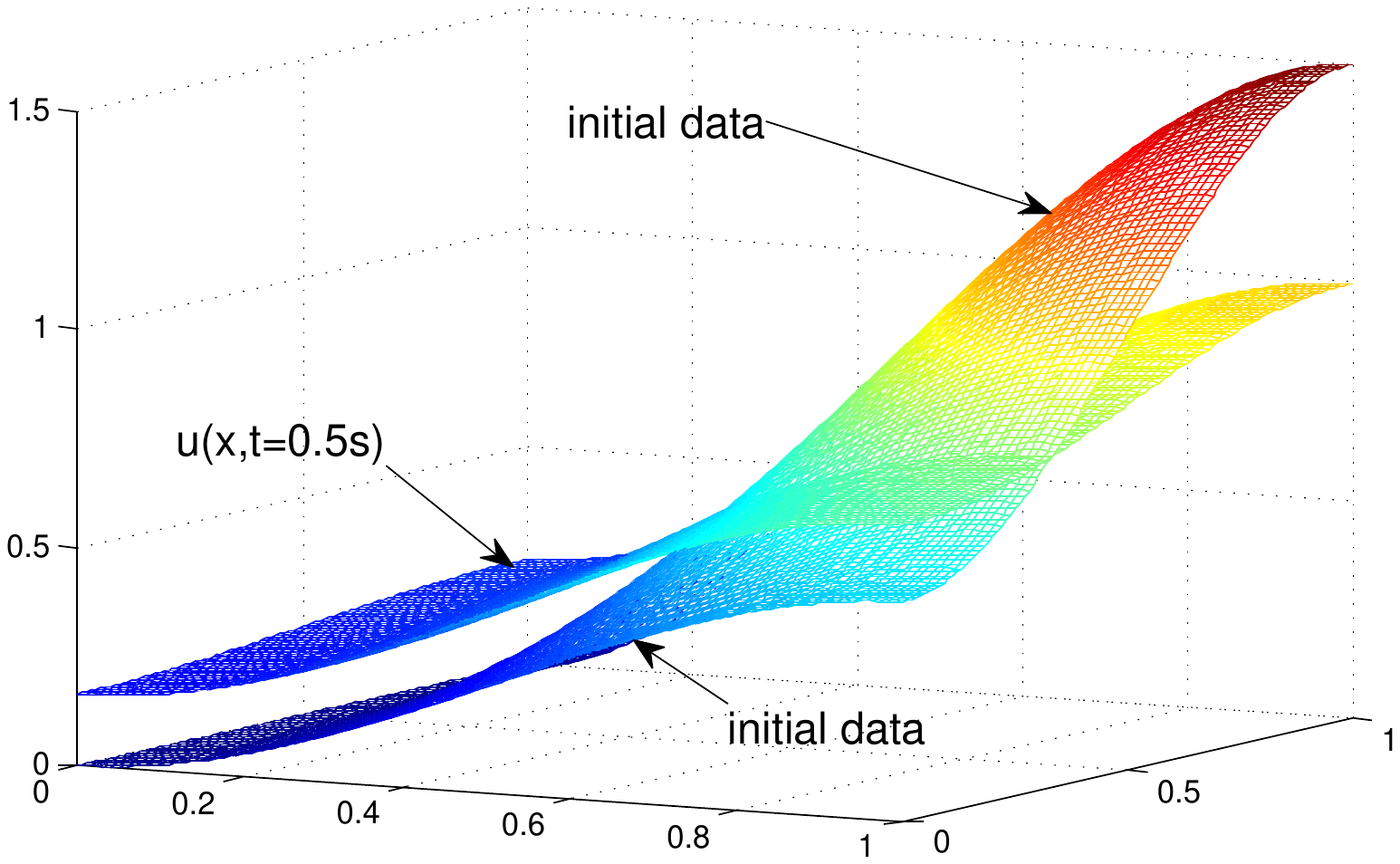}
} \hspace{10pt} \subfloat[$u(x,y)$ at t=2s] {
     \includegraphics[height=7cm,width=7cm]{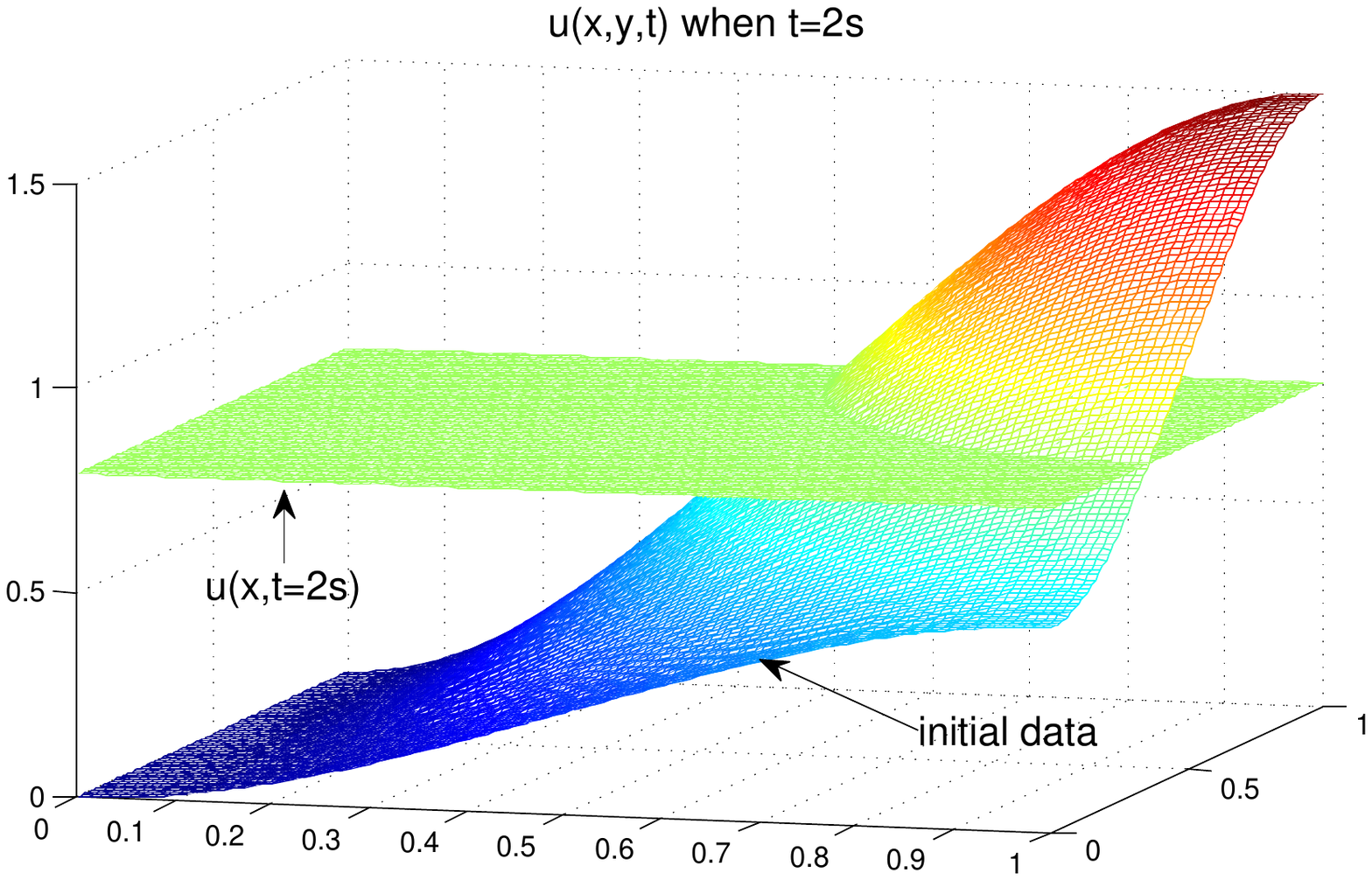}
} \hspace{10pt} \subfloat[$u(x,y)$ at t=7s] {
     \includegraphics[height=7cm,width=7cm]{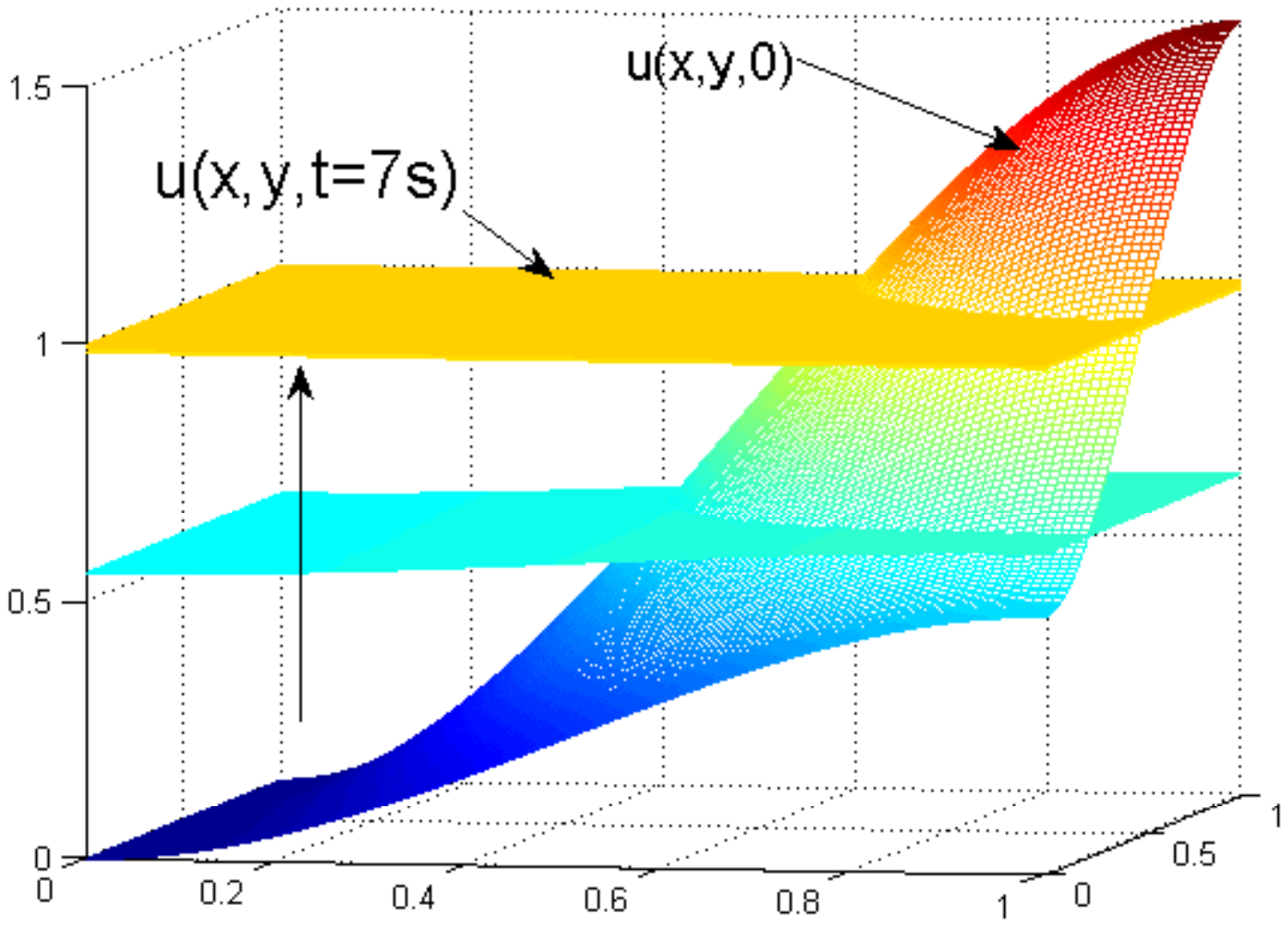}
} \hspace{10pt} \subfloat[max of $u(x,y)$ with time evolution] {
     \includegraphics[height=7cm,width=7cm]{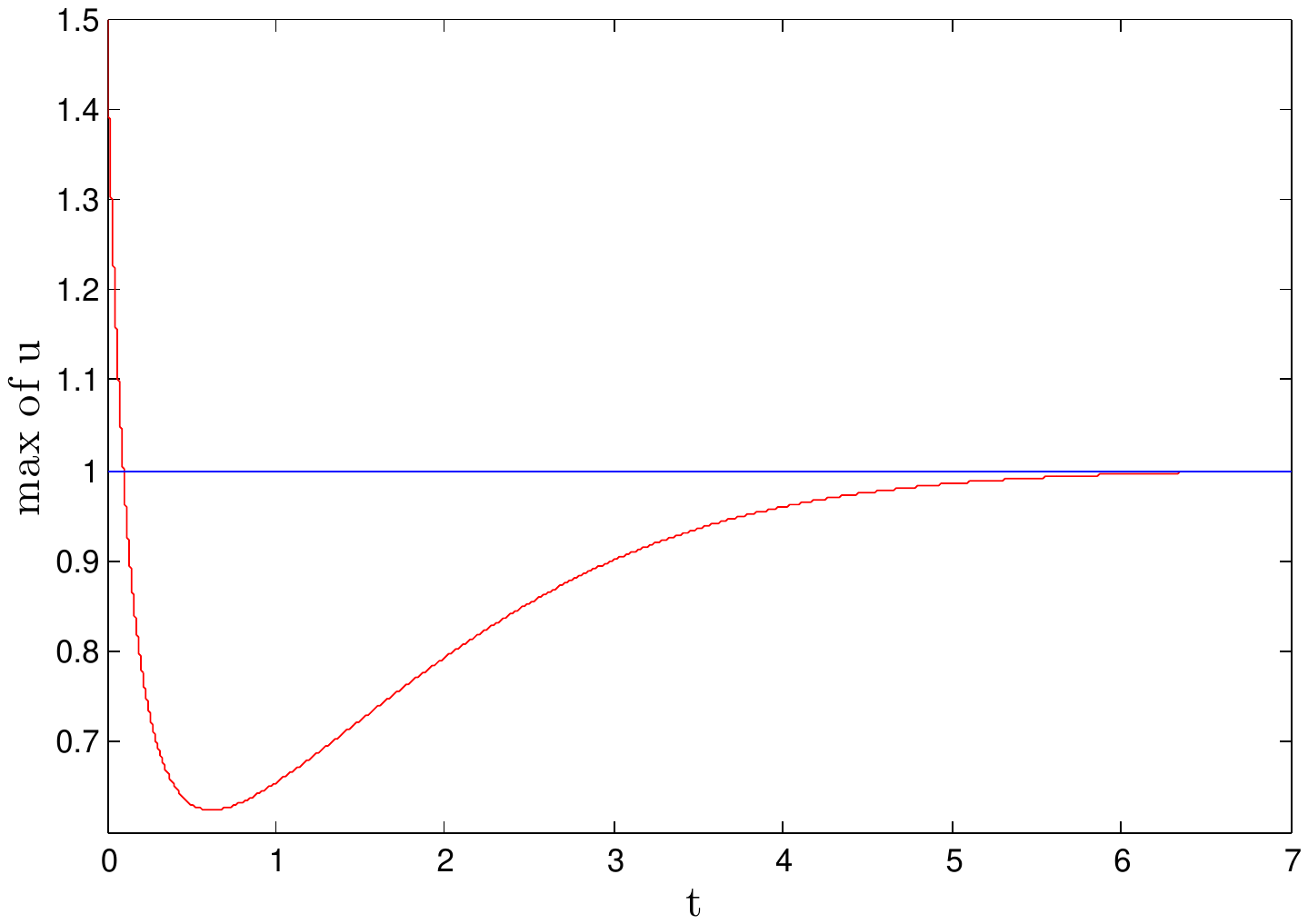}
} \caption {$u(x,y)$ with time evolution, initial mass $m_0<1$}
\label{le1}
\end{figure}

\begin{figure}[htbp]
\centering \subfloat[mass of $u(x,y)$ with time evolution] {
\includegraphics[height=6cm,width=10cm]{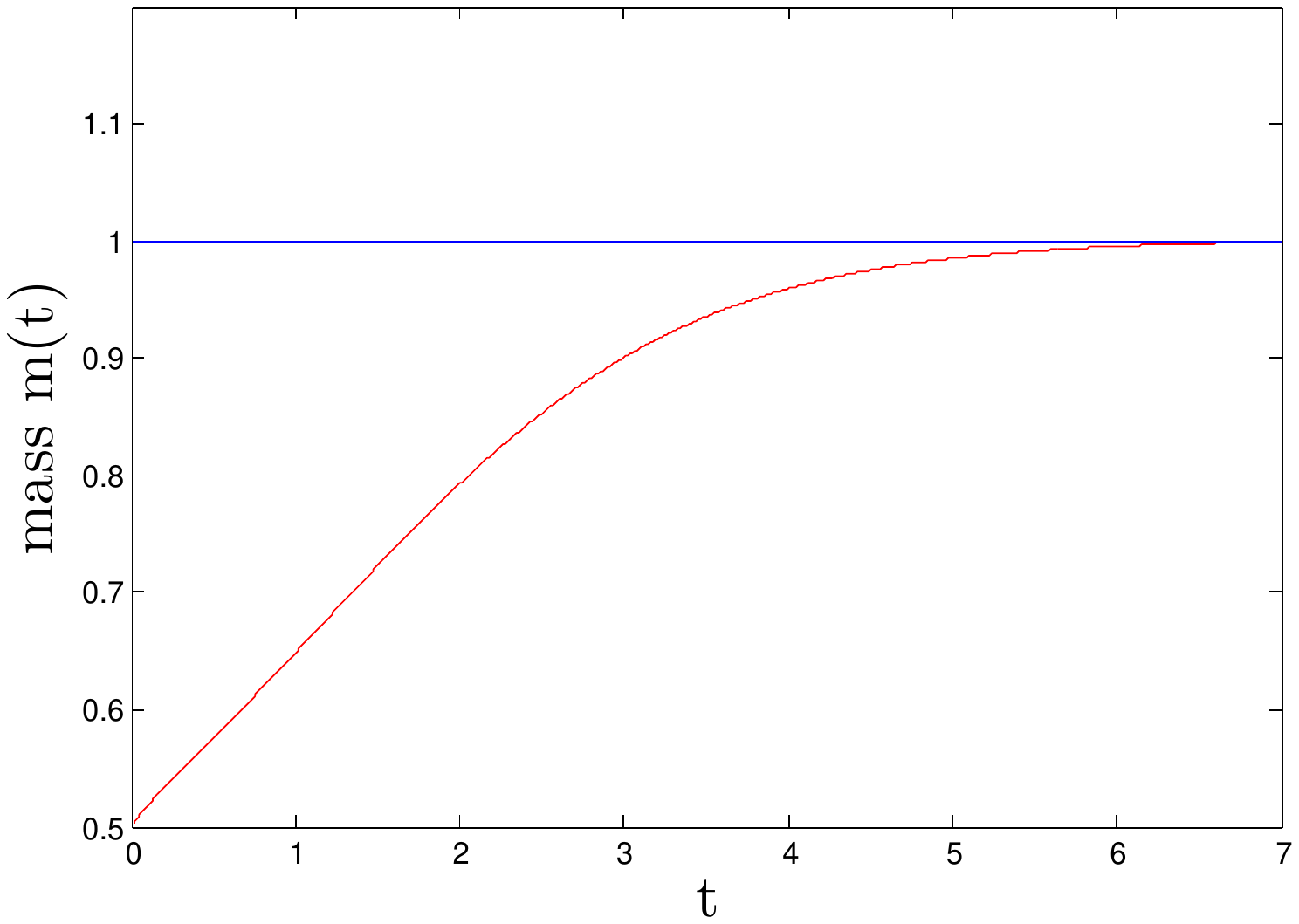} } \caption
{Mass $u$ with time evolution, initial mass $m_0<1$ }
\label{le1mass}
\end{figure}

\begin{figure}[htbp] \centering
 \subfloat[$u$ with time evolution] {
     \includegraphics[height=8cm,width=8cm]{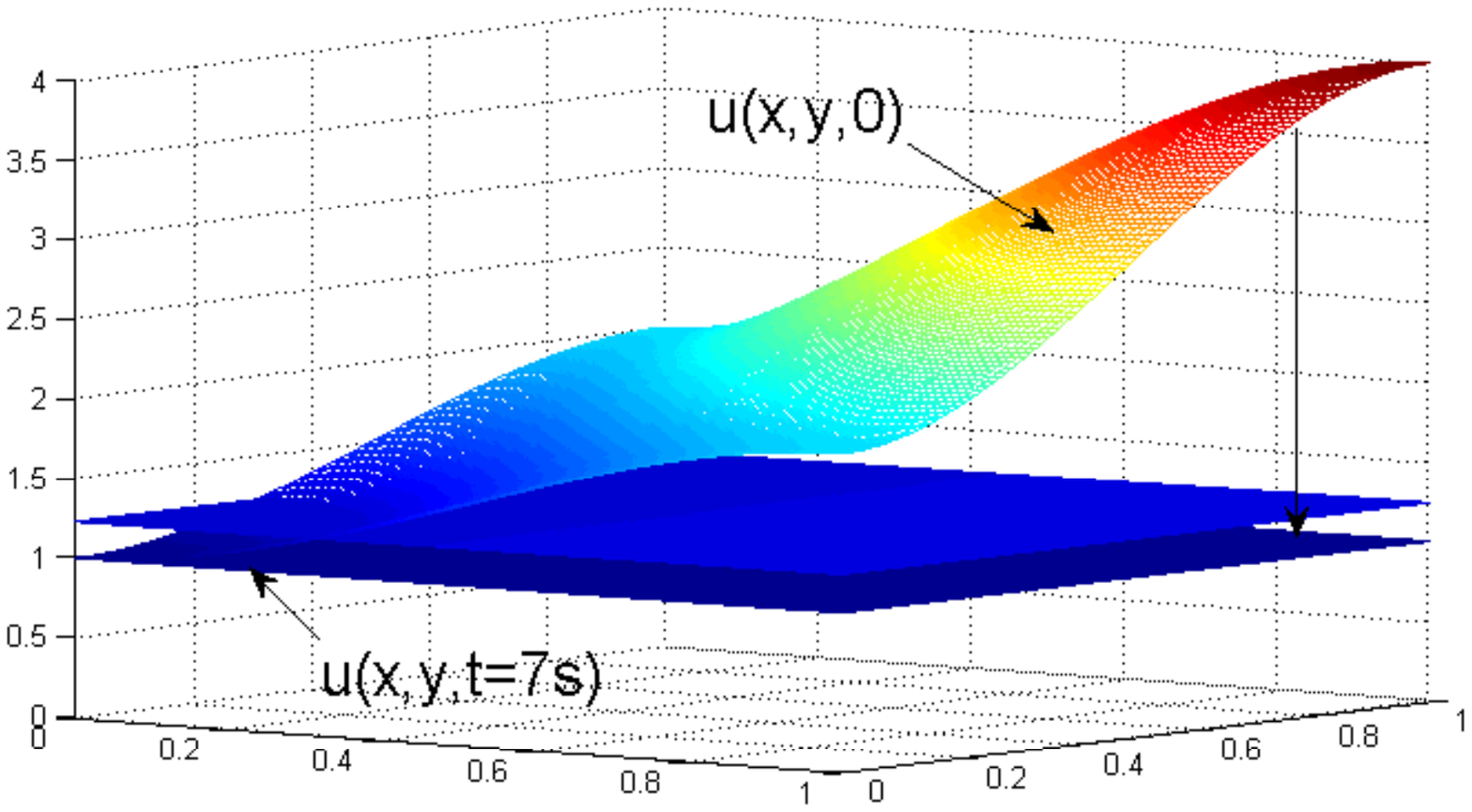}
} \hspace{10pt} \subfloat[max of $u$ with time evolution] {
     \includegraphics[height=8cm,width=8cm]{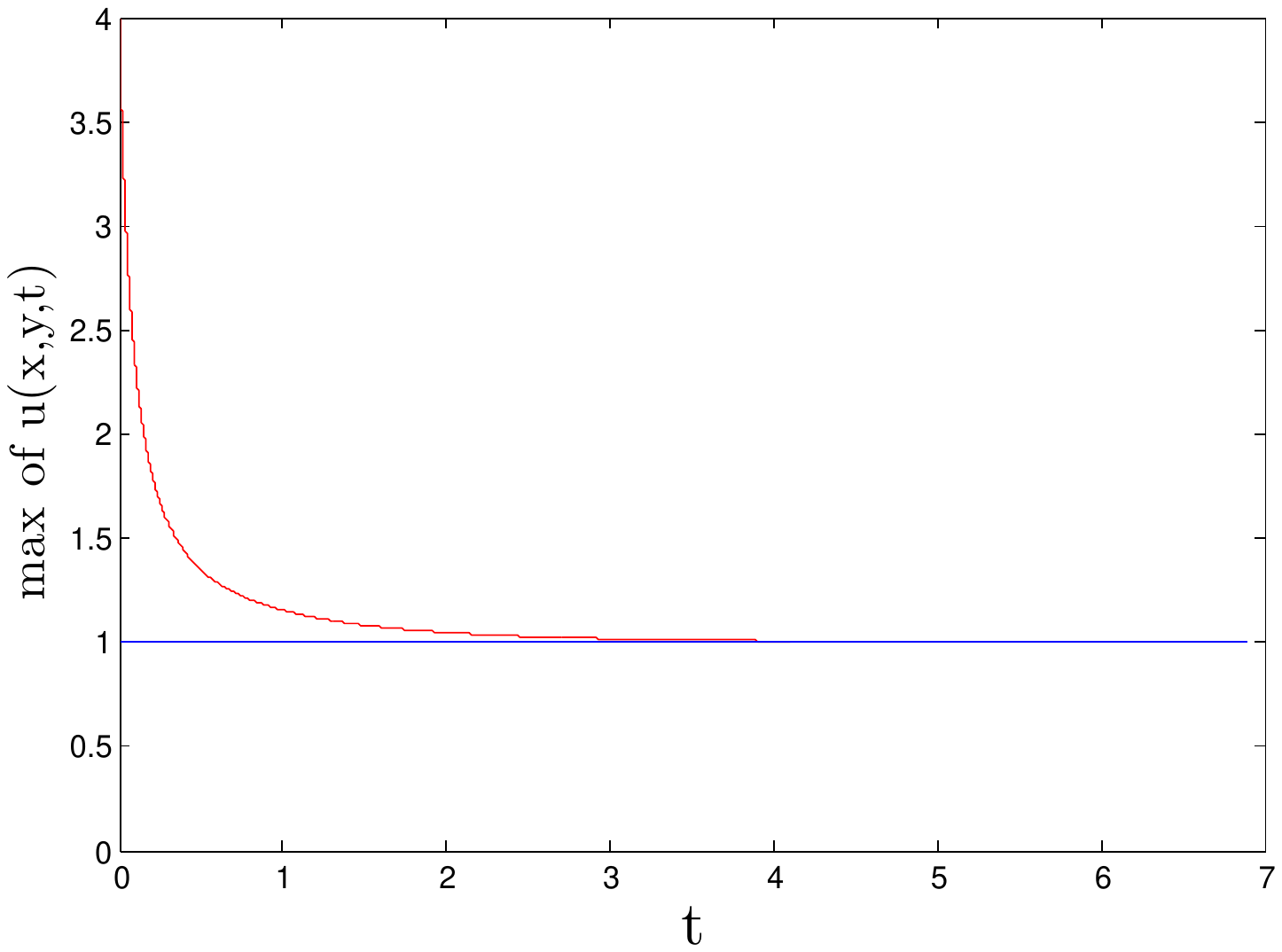}
}  \caption {$u(x,y)$ with time evolution, initial mass $m_0>1$}
\label{ge1}
\end{figure}

\begin{figure}[htbp]
\centering \subfloat[$\iint_{\Omega}u(x,y)dxdy$ with time evolution]
{ \includegraphics[height=6cm,width=10cm]{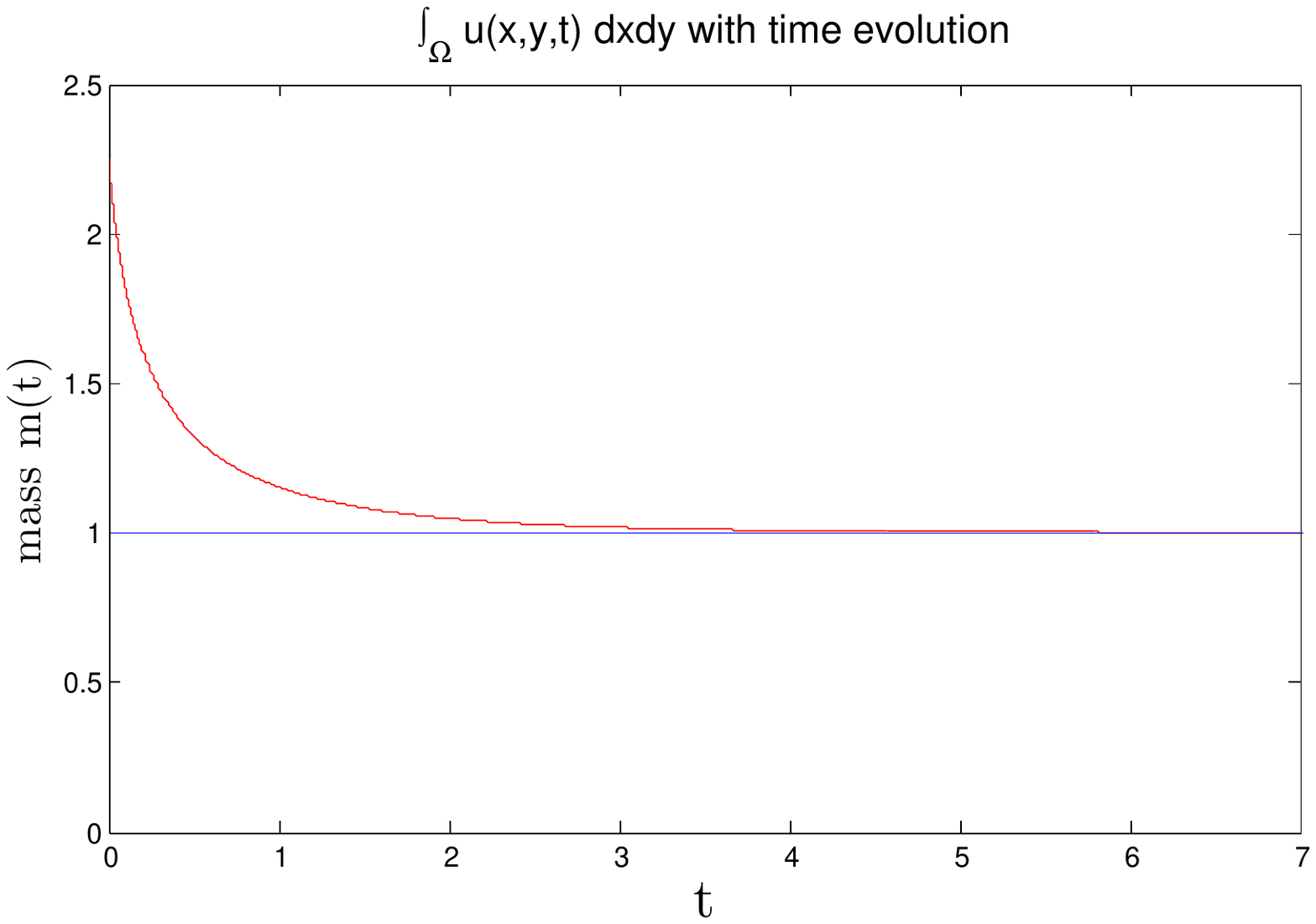} } \caption
{Mass of $u(x,y)$ with time evolution, initial mass $m_0>1$ }
\label{ge1mass}
\end{figure}

\begin{figure}[htbp] \centering
 \subfloat[initial $u(x,y)$] {
     \includegraphics[height=7cm,width=7cm]{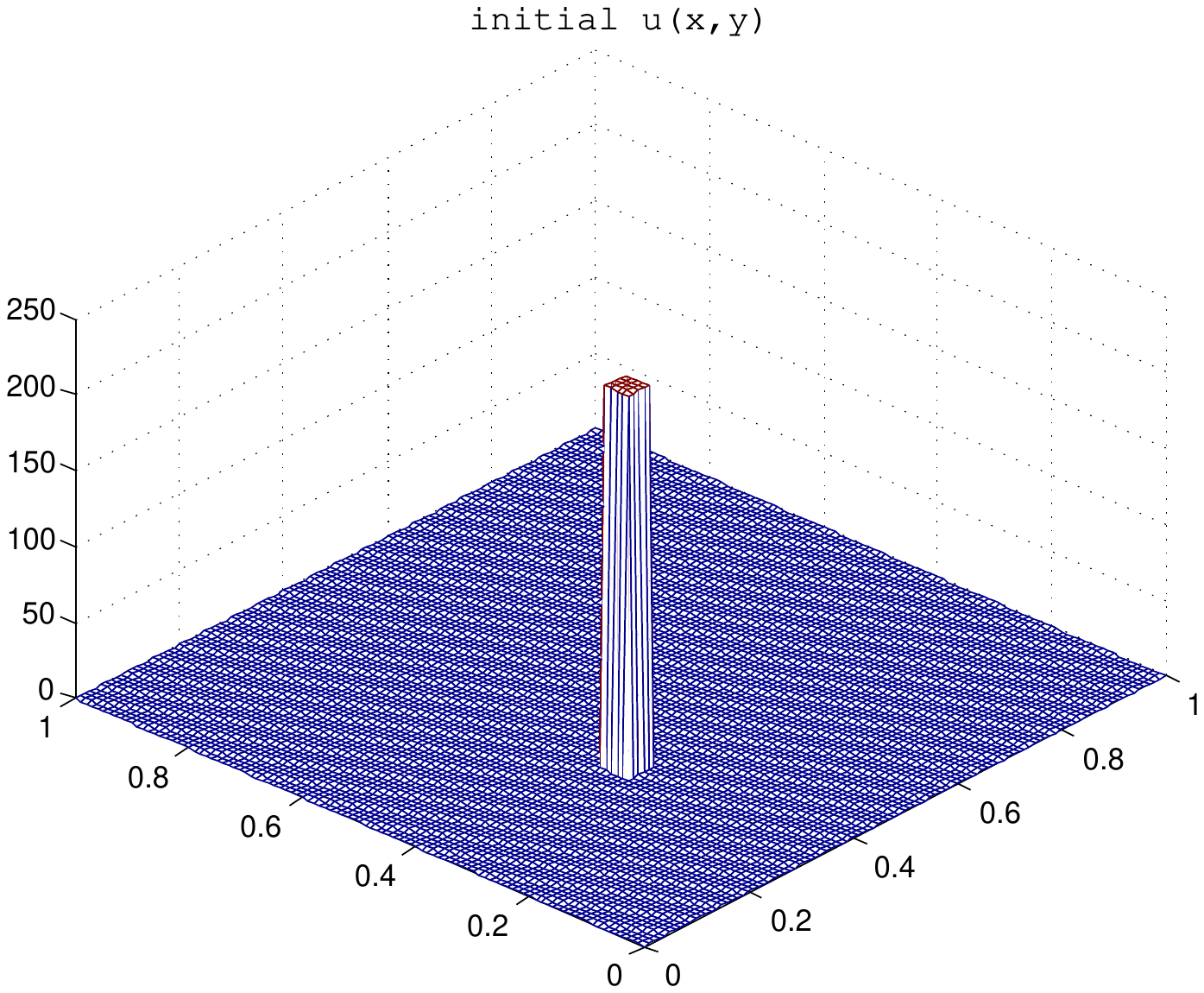}
} \hspace{10pt} \subfloat[$u(x,y)$ at $t=10^{-4}$s] {
     \includegraphics[height=7cm,width=7cm]{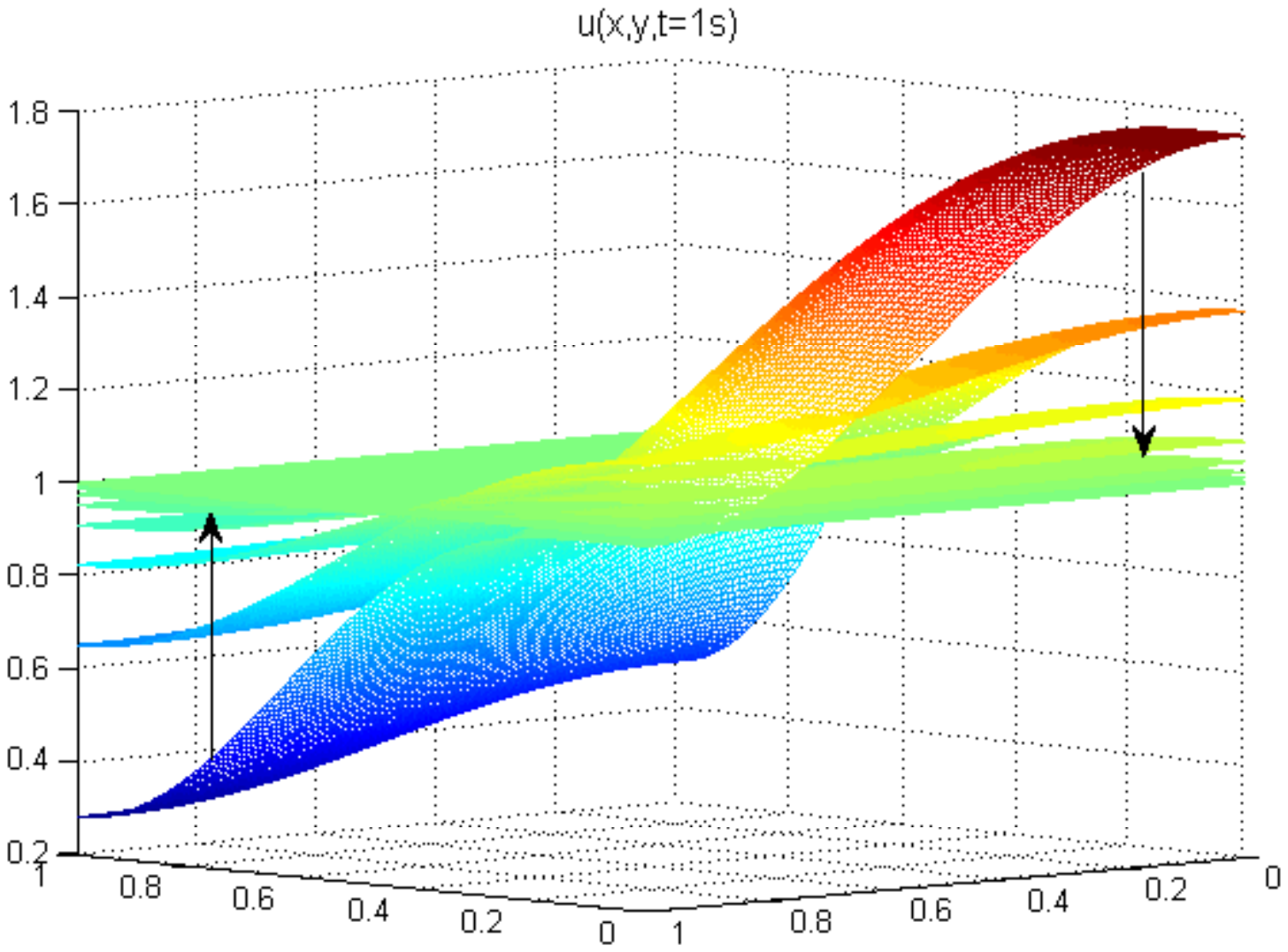}
} \hspace{10pt} \subfloat[mass of $u(x,y)$ with time evolution] {
     \includegraphics[height=7cm,width=7cm]{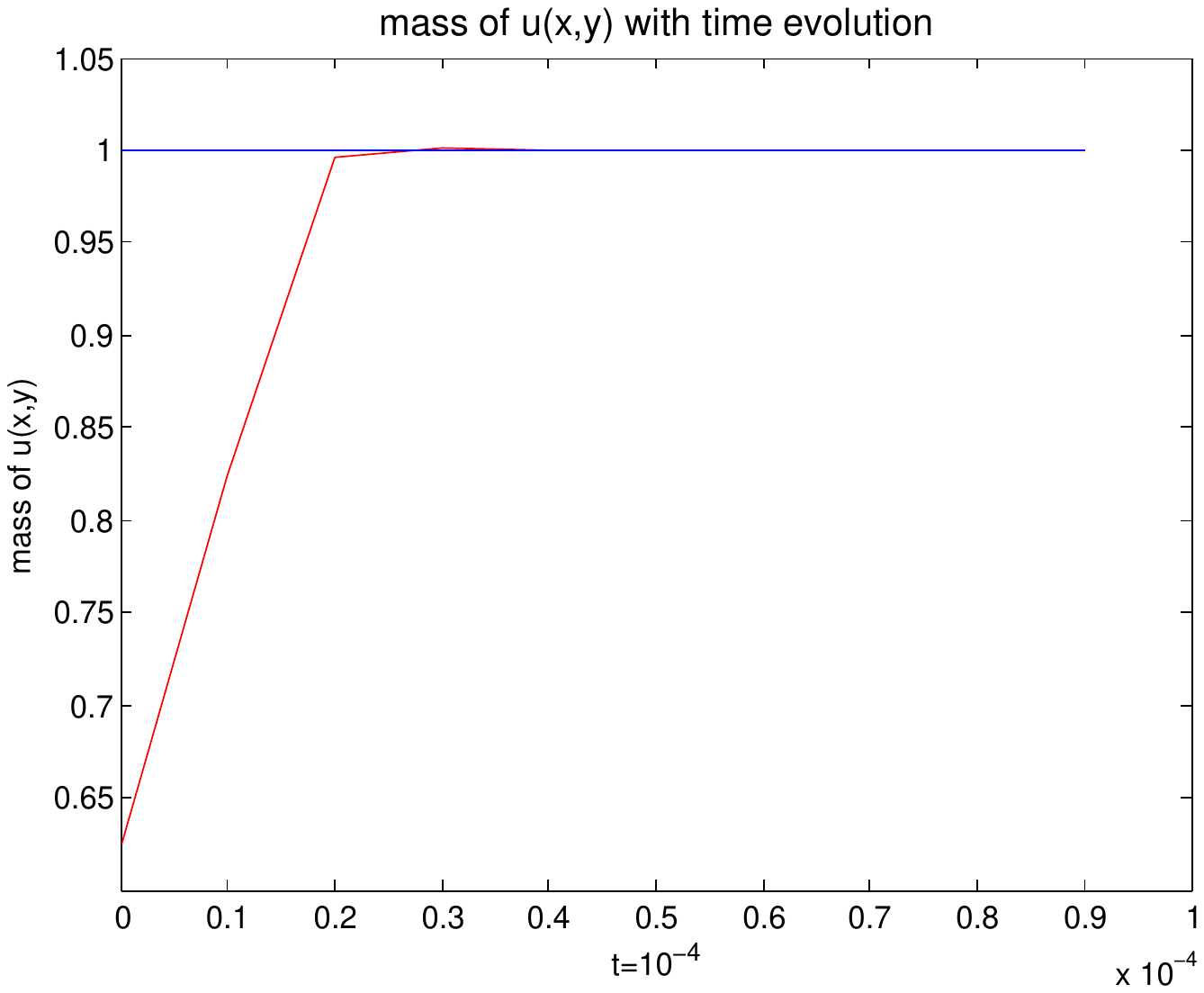}
} \hspace{10pt} \subfloat[max of $u(x,y)$ with time evolution] {
     \includegraphics[height=7cm,width=7cm]{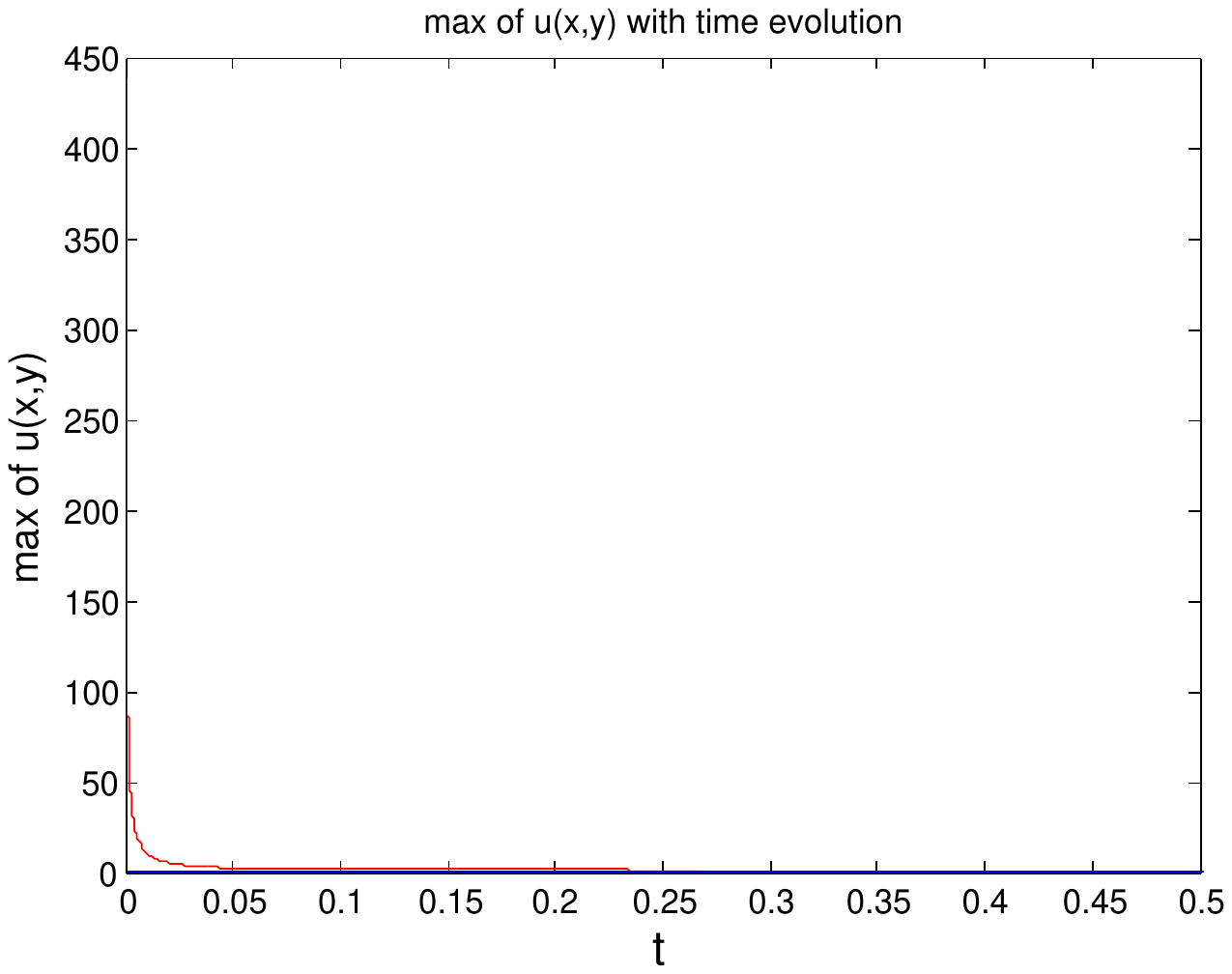}
} \caption {$n=2$, $u(x,y)$ with time evolution, $m_0<1$}
\label{less1large}
\end{figure}

\begin{figure}[htbp] \centering
 \subfloat[$u(x)$ with time evolution] {
     \includegraphics[height=7cm,width=7cm]{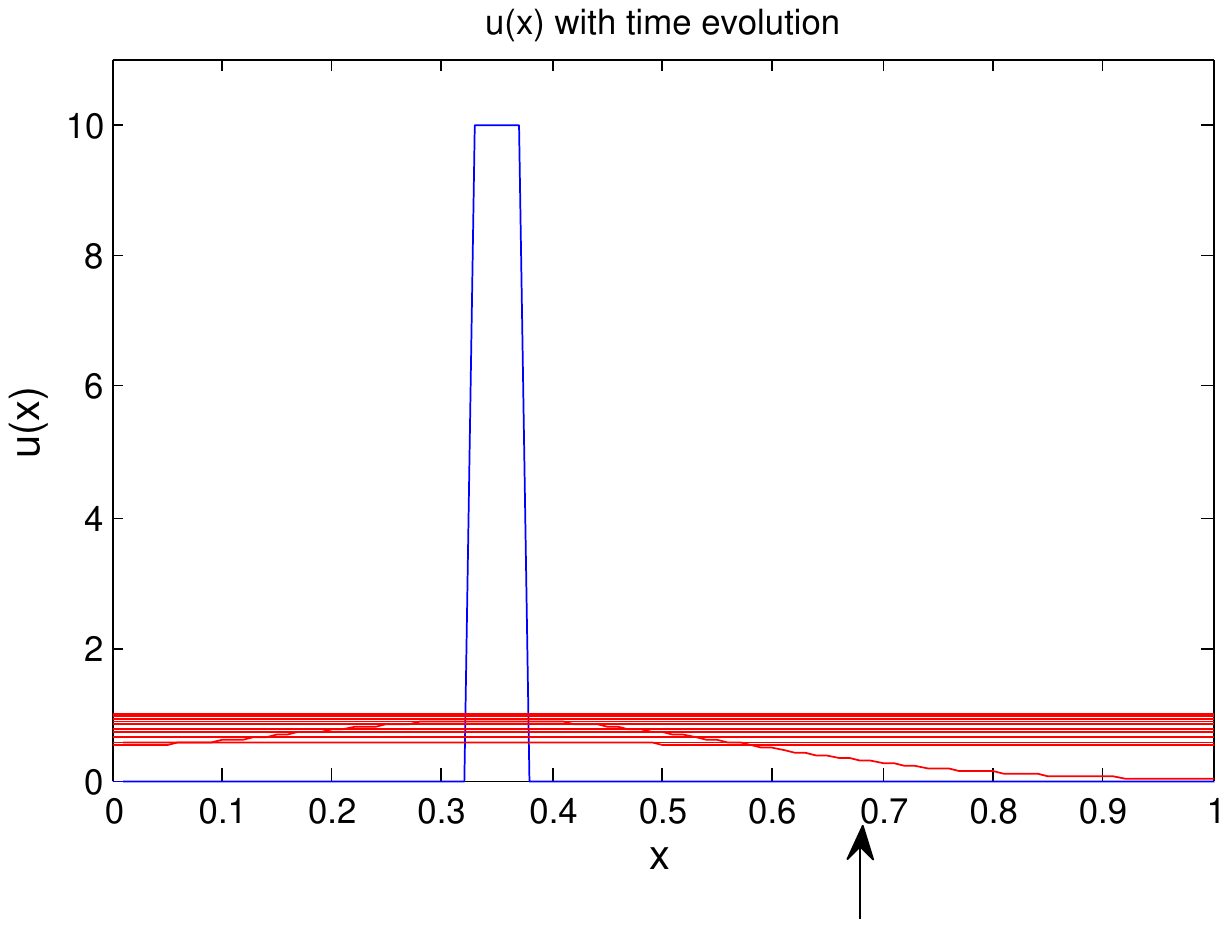}
} \hspace{10pt} \subfloat[max of $u(x)$ with time evolution] {
     \includegraphics[height=7cm,width=7cm]{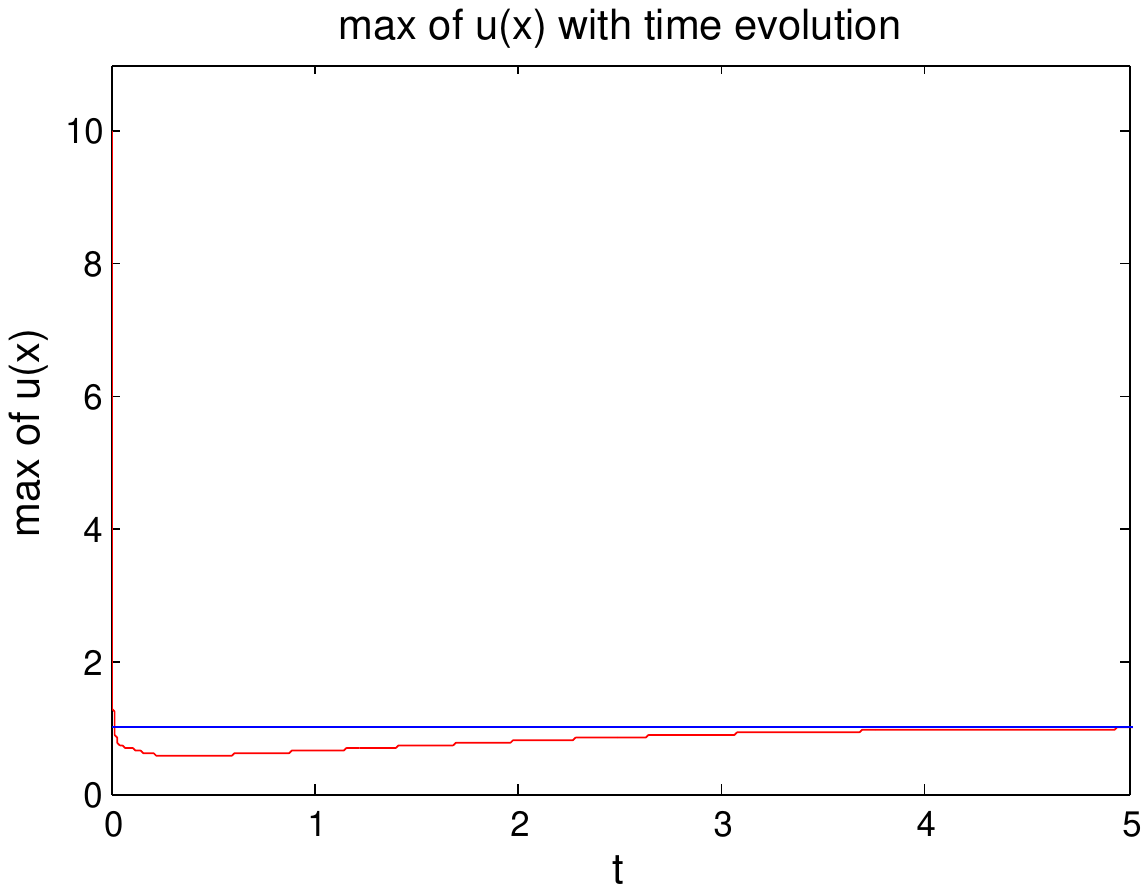}
}  \caption {$n=1$, $u(x)$ with time evolution, $m_0<1$} \label{1d}
\end{figure}

\begin{figure}[htbp]
\centering \subfloat[$\int_{\Omega}u(x)dx$ with time evolution] {
\includegraphics[height=6cm,width=10cm]{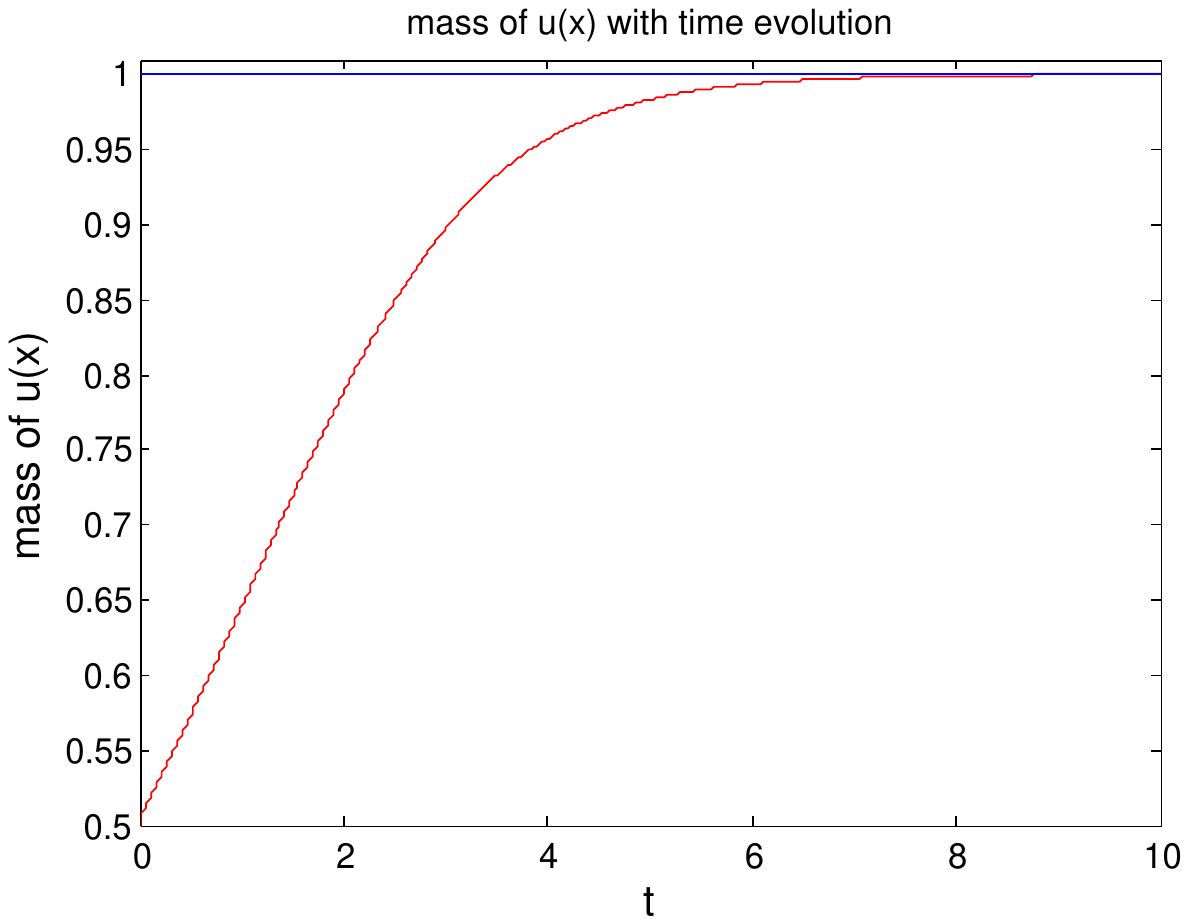} } \caption
{$n=1$, mass of $u(x)$ with time evolution, initial mass $m_0<1$ }
\label{1dle1mass}
\end{figure}

We construct the ADI scheme with operator splitting method and using
Taylor expansion linearizes $f(u)$ to approximate this semi-linear
equation as follows
\begin{align}
&\frac{u_{ij}^{k+1}-u_{ij}^k}{\tau} =  \Big( \delta_x^2+\delta_y^2
\Big) \left( \frac{u_{ij}^{k+1}+ u_{ij}^{k} }{2}\right) -
\frac{\tau}{4} \left(\delta_x^2+1 \right) \delta_y^2 \left(
u_{ij}^{k+1}- u_{ij}^{k} \right) \nonumber\\
& \quad \quad \quad \quad ~~ + f \left( u_{ij}^k \right)+
\frac{u_{ij}^{k+1}-u_{ij}^k}{2} f'\left(u_{ij}^k \right),~~,2 \le
i,j \le N-1, \label{chafen} \\
& u^1_{i,j}=U_0 \Big( x(i),y(j) \Big), ~~1 \le i,j \le N. \nonumber
\end{align}
(\ref{chafen}) can be rearranged into
\begin{align}
&\Big[  1-\frac{\tau}{2}\delta_x^2-\frac{\tau}{2}f'\left( u_{ij}^{k}
\right)    \Big]\Big[ 1- \frac{\tau}{2}\delta_y^2 \Big]
u_{ij}^{k+1}=h \left( u_{ij}^k \right), \label{hij} \\
&h \left(u_{ij}^k \right)= \Big[   1- \frac{\tau}{2}f'\left(
u_{ij}^{k}\right) \Big] u_{ij}^k + \tau f\left(u_{ij}^k\right)
+\frac{\tau}{2} \delta_x^2 u_{ij}^k + \left[  \frac{\tau^2}{4}
\delta_x^2+ \frac{\tau^2}{4} f'\left( u_{ij}^{k} \right) +
\frac{\tau}{2} \right] \delta_y^2 u_{ij}^k.  \nonumber
\end{align}
(\ref{hij}) is equivalent to the following ADI scheme
\begin{align}
&\Big[  1-\frac{\tau}{2}\delta_x^2-\frac{\tau}{2}f'\left( u_{ij}^{k}
\right)  \Big] \bar{u}_{ij}=h\left(u_{ij}^k \right),~~2 \le i \le
N-1,~
j=2,3,\cdots,N-1,~k=1,2,\cdots,K. \\
&\bar{u}_{1,j}=\frac{4 \bar{u}_{2,j}-\bar{u}_{3,j}}{3},
~~\bar{u}_{N,j}=\frac{4
\bar{u}_{N-1,j}-\bar{u}_{N-2,j}}{3}, ~~j=2,3,\cdots,N-1, \label{x0N} \\
&\Big( 1- \frac{\tau}{2}\delta_y^2 \Big)u_{ij}^{k+1} = \bar{u},~~2
\le j \le N-1,~i=2,3,\cdots,N-1,~k=1,2,\cdots,K-1, \\
&u_{i,1}^{k+1}=\frac{4 u_{i,2}^{k+1}-u_{i,3}^{k+1}}{3},~~
u_{i,N}^{k+1}=\frac{4 u_{i,N-1}^{k+1}-u_{i,N-2}^{k+1}}{3},
~~i=2,3,\cdots,N-1. \label{y0N}
\end{align}

\subsection{Numerical examples}
\indent Let $\Omega=[0,1]^2$. In the numerical study, we will consider three cases:
\begin{enumerate}
  \item $n=2, 1<\alpha<2$, we choose $m_0<1$ and $u_0 \in L^{\alpha}(\Omega)$.
  \item $n=2, \alpha>2$, we choose $m_0<1$ and $u_0 \in L^{\alpha}(\Omega)$, $\int_{\Omega}u_0^{\alpha} dx$ is large enough.
  \item $n=1,\alpha=2$, we choose $m_0<1$ and $u_0 \in L^{\alpha}(\Omega)$, $\int_{\Omega}u_0^{2} dx$ is large enough.
\end{enumerate}

\textbf{Case 1:} $n=2$, we choose $\alpha=3/2$ and the initial data
\begin{align}
u_0(x,y)=\left(-2x^3+3x^2+0.5\right)\left(-2y^3+3y^2\right),~~\iint_{\Omega}
u_0(x,y) dxdy=0.5<1.
\end{align}
Fig. \ref{le1} shows the evolution of the solutions with time. We
notice that the solution converges to $1/|\Omega|$, where $|\Omega|$
is the area of the 2-n domain. Fig. \ref{le1mass} shows its
corresponding total mass with time evolution. Eventually the mass
converges to 1.

For the initial mass is greater than 1, we choose
\begin{align}
u_0(x,y)=\left(-2x^3+3x^2+1\right)\left(-2y^3+3y^2+1\right),~~\iint_{\Omega}
u_0(x,y) dxdy=2.25>1.
\end{align}
the results are shown in Fig. \ref{ge1}, where the solution
converges to $1/|\Omega|$. The evolution of mass is shown in Fig.
\ref{ge1mass}, we observe that it finally converges to 1.

\textbf{Case 2:} $n=2$, $\alpha=3$ and $m_0<1$, the initial data is chosen to be a characteristic function
\begin{align*}
u=\left\{
    \begin{array}{ll}
     \frac{1}{40 h^2} , & 0.3 \le x \le 0.3+5h,~ 0.3 \le y \le 0.3+5h, \\
      0, & \mbox{other},
    \end{array}
  \right.
\end{align*}
where $h=0.01$. Simple computations deduce $\int_{\Omega}u_0 dx =m_0=0.625<1$ and $\int_{\Omega} u_0^\alpha dx=3.9 \times 10^{4}$. Fig. \ref{less1large} shows
the evolution of the solutions with time. We can notice that the mass tends to one very quickly and $u(x,y)$ goes to $1/|\Omega|$ with time evolution.

\textbf{Case 3:} $n=1$, $\alpha=2$ and $m_0<1$, the initial data is chosen to be
\begin{align}
u=\left\{
    \begin{array}{ll}
     10 , & 0.3 \le x \le 0.35, \\
      0, & \mbox{other},
    \end{array}
  \right.
\end{align}
the initial mass $m_0=0.5$. It can be observed in Fig. \ref{1d} that $u(x)$ converges to the constant $1/|\Omega|$, Fig. \ref{1dle1mass} shows the mass tends to one finally.

\section{Conclusions}
This paper concerns the nonlocal Fisher KPP problem (\ref{nkpps0}) with the reaction term's power $\alpha$. For $\alpha \ge 1$, the global existence of a classical solution to (\ref{nkpps0}) is analyzed. When the initial mass is less than one, $1 \le \alpha<2$ for $n=1,2$ or $1\le \alpha <1+2/n$ for $n \ge 3$, there exists a global unique nonnegative classical solution. When the initial mass is greater than or equal to one, the Fisher KPP problem admits a unique classical solution for any $ \alpha \ge 1$. Our numerical simulations show that when the initial mass is less than one and $\alpha >2$ for $n=1,2$, the unique nonnegative classical solution will exist globally. Therefore, for $n=1,2$ with $\alpha \ge 2$ or $n \ge 3$ with $\alpha \ge  2/n$, our conjecture is that the problem (\ref{nkpps0}) also admits a global unique nonnegative classical solution. This is a challenging problem which we will study in the future.



\begin{thebibliography}{9}


\bibitem{akl10}M. Anguiano, P.E. Kloeden, T. Lorenz. Asymptotic behaviour of nonlocal reaction-diffusion equations. Nonlinear Analysis, 73 (2010), 3044-3057.

\bibitem{ac94} S. N. Antontsev and M. Chipot, The thermistor problem: Existence, smoothness, unique- ness, blow-up, SIAM J. Math. Anal., 25 (1994), pp. 1128-1156.

\bibitem{ac97} S. N. Antontsev and M. Chipot, The analysis of blow-up for the thermistor problem, Siberian Math. J., 38 (1997), pp. 827-841.

\bibitem{bal1} J. M. Ball, Remarks on blow--up and nonexistence theorems for
nonlinear evolution equations, Quart. J. Math. Oxford, (2), {\bf 28}
(1977), 473--486.

\bibitem{bal2} J. M. Ball, Finite time blow--up in nonlinear problems,
Nonlinear Evolution Equations (1977), Academic Press, 189--205.

\bibitem{bb82} J. W. Bebernes and A. Bressan, Thermal behaviour for a confined reactive gas, J. Differential Equations, 44 (1982), pp. 118-133.

\bibitem{be82} J. W. Bebernes and R. Ely, Comparison techniques and the method of lines for a parabolic functional equation, Rocky Mountain J. Math., 12 (1982), pp. 723-733.

\bibitem{bt96}J. W. Bebernes and P. Talaga, Nonlocal problems modelling shear banding, Comm. Appl. Nonlinear Anal., 3 (1996), pp. 79-103.

\bibitem{BL13}
\newblock S. Bian and J.-G. Liu,
\newblock \emph{Dynamic and steady states for multi-dimensional Keller-Segel model with diffusion exponent $m > 0$},
\newblock Comm Math Phy., \textbf{323} (2013), 1017-1070.


\bibitem{BL14}
\newblock S. Bian, J.-G. Liu and C. Zou,
\newblock Ultra-contractivity for Keller-Segel model with diffusion exponent $m>1-2/d$,
\newblock Kinetic and Related Models, \textbf{7}(1) (2014), 9-28.

\bibitem{bds93}C. Budd, B. Dold, and A. Stewart, Blowup in a partial differential equation with con- served first integral, SIAM J. Appl. Math., 53 (1993), pp. 718-742.

\bibitem{d95}  K. Deng, Dynamical behavior of solutions of a semilinear parabolic equation with nonlocal
singularity, SIAM J. Math. Anal., 26 (1995), pp. 98-111.

\bibitem{dkl92} K. Deng, M. K. Kwang, and H. A. Levine, The influence of nonlocal nonlinearities on the
long time behaviour of solutions of Burgers equation, Quart. Appl. Math., 50 (1992),
pp. 173-200.

\bibitem{dlx03}W. Deng, Y. Li, C. Xie. Semilinear reaction-diffusion systems with nonlocal sources, Mathematical and Computer Modelling, 37 (2003), 937-943.

\bibitem{f97}M. Fila, Boundedness of global solutions of nonlocal parabolic equations, Nonlinear Anal.,
30 (1997), pp. 877-885.

\bibitem{f37}Fisher RA. The wave of advance of advantageous genes. Ann Eugenics 1937, (7), 355-69.

\bibitem{F69}
\newblock A. Friedman,
\newblock ``The Porous Medium Equation: Mathematical Theory'',
\newblock Oxford University Press, 2007.


\bibitem{fuj1} H. Fujita, On the nonlinear equations $\Delta u+e^u=0$ and $v_t=\Delta v+e^v$, Bull. Amer. Math. Soc., 75 (1969), 132--135.

\bibitem{fuj2} H. Fujita, On some nonexistence and nonuniqueness theorems
for nonlinear parabolic equations, Proc. Symp. Pure Math. XVIII,
Nonlinear Functional Analysis Amer. Math. Soc., {\bf 28} (1970),
105--113.

\bibitem{fuj3} H. Fujita, On the blowing up of solutions of the Cauchy problem for $u_t=\Delta u+u^{1+\alpha}$, J. Fac. Sci. Univ. Tokyo Sect. IA Math. {\bf13} (1966) 109--124.

\bibitem{SimJac86}Simon, Jacques. "Compact sets in the spaceL p (O, T; B)." Annali di Matematica pura ed applicata 146.1 (1986): 65-96.

\bibitem{hy95}B. Hu and H.Ã¢ÂÂM. Yin, Semilinear parabolic equations with prescribed energy, Rend. Circ.
Mat. Palermo, 44 (1995), pp. 479-505.

\bibitem{kapl} S. Kaplan, On the growth of solutions of quasilinear
parabolic equations, Comm. Pure Appl. Math. {\bf 16} (1963), 305--330.

\bibitem{kpp}Kolmogorov AN, Petrovsky IG, Piskunov NS., Investigation of the equation of diffusion combined with increasing of the substance and its
application to a biology problem. Bull Moscow State Univ Ser A: Math and Mech 1937;1(6):1-25.

\bibitem{l95ab}A. A. Lacey, Thermal runaway in a non-local problem modelling ohmic heating. I: Model
derivation and some special cases, European J. Appl. Math., 6 (1995), pp. 127Ã¢ÂÂ144; II: General proofs of blow-up and asymptotics of runaway, European J. Appl. Math., 6 (1995), pp. 201-224.

\bibitem{lsv} O. A. Ladyzenskaja, V. A. Solonnikov and N. N. Ural'ceva,
Linear and Quasilinear Equations of Parabolic Type, Transl. Math.
Monog. Vol. 23, Amer. Math. Soc., Providence, R.I., 1968.

\bibitem{lb96}G.M. Lieberman, {\it Second Order Parabolic Partial Differential Equations}, World Scientific, 1996.

\bibitem{lcl09}Q. Liu, Y. Chen, S. Lu. Uniform blow-up profiles for nonlinear and nonlocal reaction-diffusion equations. Nonlinear Analysis, 71 (2009), 1572-1583.

\bibitem{Lorz:2011hl} Lorz, Alexander and Mirrahimi, Sepideh and Perthame, Benoit, Dirac mass dynamics in multidimensional nonlocal parabolic equations, Commun Part Diff Eq, 2011, vol.36, 6, pa.1071--1098.

\bibitem{Lorz:2013vp} Lorz, Alexander and Lorenzi, Tommaso and Clairambault, Jean and Escargueil, Alexandre and Perthame, Benoit, Effects of space structure and combination therapies on phenotypic heterogeneity and drug resistance in solid tumors, arXiv:1312.6237v1 [q-bio.TO]


\bibitem{Lorz:2013hq} Lorz, Alexander and Lorenzi, Tommaso and Hochberg, Michael E and Clairambault, Jean and Perthame, Benoit, Populational adaptive evolution, chemotherapeutic resistance and multiple anti-cancer therapies, ESAIM: M2AN, 2013, vol. 47, 2, p.377--399.


\bibitem{lu95}A. Lunardi, {\it Analytic Semigroups and Optimal Regularity in Parabolic Problems}, Progr. Nonlinear Differential Equations Appl., Birkh\"auser-
Verlag, 1995.

\bibitem{Pao92}C. V. Pao, Blowing-up of solution for a nonlocal reaction-diffusion problem in combustion theory, J. Math. Anal. Appl., 166 (1992), pp. 591-600.

\bibitem{pazy}A. Pazy, Semigroups of Linear Operators and Applications to Partial Differential Equations, Appl. Math. Sci., vol. 44, Springer-Verlag, 1983.

\bibitem{sbook}Pavol Quittner, Philippe Souplet, Superlinear parabolic problems. Blow-up, global existence and steady states,
Birkh\"auser Advanced Texts, 2007, 584 p.+xi. ISBN: 978-3-7643-8441-8.


\bibitem{s98}P. Souplet, Blow-up in nonlocal reaction-diffusion equations, Siam J. Math. Anal. Vol. 29, No. 6, pp. 1301-1334, November 1998.

\bibitem{r03}P. Rouchon, Universal bounds for global solutions of a diffusion equation with a nonlocal reaction term, J. Differential Equations, 193 (2003), 75-94.

\bibitem{st80}H.B. Stewart, Generation of analytic semigroups by strongly elliptic operators under general boundary conditions, Trans. Amer. Math.
Soc. 259 (1) (1980) 299-310.

\bibitem{vol1} V. Volpert. Elliptic partial differential equations. Volume 1. Fredholm theory of elliptic problems in unbounded domains. Birkh\"auser, 2011.
\bibitem{vol2} V. Volpert. Elliptic partial differential equations. Volume 2. Reaction-diffusion equations. Birkh\"auser, 2014.


\bibitem{volpet} V. Volpert, S. Petrovskii, Reaction-diffusion waves in biology, Physics of Life Reviews 6 (2009) 267-310.

\bibitem{VVpp}V. Volpert, V. Vougalter, Existence of stationary pulses for nonlocal reaction-diffusion equations, preprint.

\bibitem{XiuAnsgar}Chen, Xiuqing, Ansgar JÃÂ¼ngel, and Jian-Guo Liu. "A Note on Aubin-Lions-Dubinski ÃÂ­ Lemmas." Acta Applicandae Mathematicae (2013): 1-11.

\bibitem{ww11}X. Wang, W. Wo. Long time behavior of solutions for a scalar nonlocal reaction-diffusion equation. Arch. Math., 96 (2011), 483-490.

\bibitem{ww96} M. Wang and Y. Wang, Properties of positive solutions for non-local reaction-diffusion problems, Math. Methods Appl. Sci., 19 (1996), pp. 1141-1156.

\bibitem{Yi96}Y. Yin, Quenching for solutions of some parabolic equations with singular nonlocal terms, Dynam. Systems Appl., 5 (1996), pp. 19-30.

\bibitem{zfk38}Zeldovich YaB, Frank-Kamenetskii DA. A theory of thermal propagation of flame. Acta Physicochim USSR 1938;9:341-50.

\end{thebibliography}
\end{document}